\documentclass{article}

\usepackage{latexsym,amssymb,cite}
\input amssym.def \input amssym

\newtheorem{te}{Theorem}[section]

\newtheorem{de}[te]{Definition}
\newtheorem{lm}[te]{Lemma}

\newtheorem{pp}[te]{Proposition}
\newtheorem{co}[te]{Corollary}
\newtheorem{ex}[te]{Example}
\newtheorem{qu}[te]{Question}

\def\dokaz{\noindent{\bf Proof. }}
\def\kraj{\hfill $\Box$ \par \vspace*{2mm} }
\def\widemid{\hspace{1mm}\widetilde{\mid}\hspace{1mm}}
\def\nwidemid{\hspace{1mm}\widetilde{\nmid}\hspace{1mm}}
\newcommand{\zve}[1]{{{}^*\hspace{-0.5mm}#1}}
\newcommand{\zvez}[1]{{{}^*\hspace{-1mm}#1}}
\def\zvepar{\hspace{1mm}\zvez\parallel\hspace{1mm}}
\def\nwidemid{\hspace{1mm}\widetilde{\nmid}\hspace{1mm}}
\def\zvemid{\hspace{1mm}\zvez\mid\hspace{1mm}}
\def\nzvemid{\hspace{1mm}\zvez\nmid\hspace{1mm}}

\begin{document}

\begin{center}
           {\huge \bf $\widetilde{\mid}\hspace{1mm}$-divisibility of ultrafilters II: The big picture}
\end{center}
\begin{center}
{\small \bf Boris  \v Sobot}\\[2mm]
{\small  Department of Mathematics and Informatics, University of Novi Sad,\\
Trg Dositeja Obradovi\'ca 4, 21000 Novi Sad, Serbia\\
e-mail: sobot@dmi.uns.ac.rs\\
ORCID: 0000-0002-4848-0678}
\end{center}

\begin{abstract}
A divisibility relation on ultrafilters is defined as follows: ${\cal F}\hspace{1mm}\widetilde{\mid}\hspace{1mm}{\cal G}$ if and only if every set in $\cal F$ upward closed for divisibility also belongs to $\cal G$. After describing the first $\omega$ levels of this quasiorder, in this paper we generalize the process of determining the basic divisors of an ultrafilter. First we describe these basic divisors, obtained as (equivalence classes of) powers of prime ultrafilters. Using methods of nonstandard analysis we define the pattern of an ultrafilter: the collection of its basic divisors as well as the multiplicity of each of them. All such patterns have a certain closure property in an appropriate topology. We isolate the family of sets belonging to every ultrafilter with a given pattern. We show that every pattern with the closure property is realized by an ultrafilter. Finally, we apply patterns to obtain an equivalent condition for an ultrafilter to be self-divisible.
\end{abstract}



\section{Introduction}

Let $\mathbb{N}$ denote the set of natural numbers (without zero), and $\omega=\mathbb{N}\cup\{0\}$. $\beta\mathbb{N}$ is the set of all ultrafilters on $\mathbb{N}$ and, for each $n\in\mathbb{N}$, the principal ultrafilter $\{A\subseteq\mathbb{N}:n\in A\}$ is identified with $n$. Considering the topology with base sets $\overline{A}=\{{\cal F}\in\beta\mathbb{N}:A\in{\cal F}\}$, we think of $\beta\mathbb{N}$ as an extension of the discrete space $\mathbb{N}$, called the Stone-\v Cech compactification of $\mathbb{N}$. In general, for $S\subseteq\beta\mathbb{N}$, ${\rm cl}S$ will denote the closure of $S$ in this topology; in particular, for $A\subseteq\mathbb{N}$, ${\rm cl}A=\overline{A}$.

One of the main features of $\beta\mathbb{N}$ is that each function $f:\mathbb{N}\rightarrow\mathbb{N}$ can be uniquely extended to a continuous $\widetilde{f}:\beta\mathbb{N}\rightarrow\beta\mathbb{N}$. Using this, binary operations can also be extended, so by applying this to the multiplication on $\mathbb{N}$ (and denoting the extension also by $\cdot$) a right-topological semigroup $(\beta\mathbb{N},\cdot)$ is obtained, where
\begin{equation}\label{eqtensor}
{\cal F}\cdot{\cal G}=\{A\subseteq{\mathbb{N}}:\{n\in{\mathbb{N}}:A/n\in{\cal G}\}\in{\cal F}\},
\end{equation}
and $A/n=\{a/n:(a\in A)\land(n\mid a)\}$. Many properties of this and other semigroups obtained in this way are described in the book \cite{HS}.

In \cite{So1} several ways to extend the divisibility relation $\mid$ to $\beta\mathbb{N}$ were proposed. Some of them used directly the extension of multiplication (\ref{eqtensor}) and imitated the definition of $\mid$. However, a relation obtained in a different way proved to have much nicer properties: if
$$A\hspace{-1mm}\uparrow=\{n\in {\mathbb{N}}:(\exists a\in A)a\mid n\}$$
for $A\subseteq\mathbb{N}$ and
$${\cal U}=\{A\in P({\mathbb{N}})\setminus\{\emptyset\}:A=A\hspace{-1mm}\uparrow\},$$
let ${\cal F}\hspace{1mm}\widetilde{\mid}\hspace{1mm}{\cal G}$ if for every $A\in{\cal F}$ holds $A\hspace{-1mm}\uparrow\in{\cal G}$. It turned out that another equivalent condition is more convenient in practice:
$${\cal F}\hspace{1mm}\widetilde{\mid}\hspace{1mm}{\cal G}\Leftrightarrow{\cal F}\cap{\cal U}\subseteq{\cal G}.$$
$\widemid$ is a quasiorder, and defining
$${\cal F}=_\sim{\cal G}\Leftrightarrow({\cal F}\widemid{\cal G})\land({\cal G}\widemid{\cal F}),$$
$=_\sim$ is an equivalence relation, so $\widemid$ can be thought of as a partial order on equivalence classes $[{\cal F}]_\sim\in\beta\mathbb{N}/\hspace{-0.5mm}=_\sim$. For $n\in\mathbb{N}$, $n\widemid{\cal G}$ means that ${\cal G}=n{\cal F}$ for some ${\cal F}\in\beta\mathbb{N}$; however this is not true for $n\in\beta\mathbb{N}\setminus\mathbb{N}$ in general.

It is worth mentioning that two general (so-called canonical) ultrafilter extensions of relations were described in \cite{PS}. Among other topics, their topological properties were considered. One of them (although it was defined in a different way) actually produces $\widemid$. The other one is, in a sense, trivial when applied to the divisibility relation.

We recall some notation and introduce some more. Throughout the text, for $A,B,A_i\subseteq\mathbb{N}$ and $n\in\mathbb{N}$:
\begin{itemize}
\item $\mathbb{P}$ denotes the set of primes;
\item $\mathbb{P}^{exp}=\{p^n:p\in \mathbb{P}\land n\in\mathbb{N}\}$;
\item $A^n=\{a^n:a\in A\}$;
\item $A_1A_2\dots A_n=\{a_1a_2\dots a_n:(\forall i\in\{1,2,\dots,n\})(a_i\in A_i\land(\forall j\neq i)gcd(a_i,a_j)=1)\}$;
\item $A^{(n)}=\underbrace{A\cdot A\cdot\dots\cdot A}_n$;
\item $A\hspace{-0.1cm}\downarrow=\{n\in {\mathbb{N}}:(\exists a\in A)n\mid a\}$ and
\item ${\cal V}=\{A\in P({\mathbb{N}})\setminus\{\mathbb{N}\}:A=A\hspace{-0.1cm}\downarrow\}$.
\end{itemize}

A set $A\subseteq\mathbb{N}$ is convex if, for all $x,y\in A$ and all $z\in\mathbb{N}$ such that $x\mid z$ and $z\mid y$ we have $z\in A$. Let $\cal C$ be the family of all convex subsets of $\mathbb{N}$. Clearly, since ${\cal U}\subseteq{\cal C}$, ${\cal F}\cap{\cal C}$ determines the $=_\sim$-class of an ultrafilter $\cal F$.

Various aspects of $(\beta\mathbb{N}/\hspace{-1mm}=_\sim,\widemid)$ were considered in papers \cite{So1}--\cite{So7}. In particular, the paper \cite{So3} contains many facts on the first $\omega$ levels of this order, some of which we will now briefly recapitulate. The $n$-th level $\overline{L_n}$ consists of ultrafilters that concentrate on the set $L_n=\{p_1p_2\dots p_n:p_1,p_2,\dots,p_n\in\mathbb{P}\}$ ($p_i$'s need not be distinct). Hence, $L_0=\{1\}$ and $L_1=\mathbb{P}$. The ultrafilters in $\overline{\mathbb{P}}$ are called {\it prime}. For ${\cal P}\in\overline{\mathbb{P}}$, we denote
$${\cal P}\upharpoonright\mathbb{P}=\{A\in{\cal P}:A\subseteq\mathbb{P}\}.$$
The square ${\cal P}^2$ of ${\cal P}\in\overline{\mathbb{P}}$ is generated by the sets $A^2$ for $A\in{\cal P}\upharpoonright\mathbb{P}$. However, it is not the same as ${\cal P}\cdot{\cal P}$; the latter belongs to the family of ultrafilters ``twice divisible" by ${\cal P}$ and containing not the sets $A^2$ but $A^{(2)}$ for $A\in{\cal P}\upharpoonright\mathbb{P}$. Example \ref{exdifprod} will make this distinction more natural.

Powers of primes ${\cal P}^k$ (for ${\cal P}\in\overline{\mathbb{P}}$ and $k\in\mathbb{N}$) are called {\it basic ultrafilters}. All the elements of 
$$L:=\bigcup_{n\in\omega}\overline{L_n}$$
can be fragmented into such basic ingredients. This is not a factorization: it simply means that an ultrafilter in $\overline{L_n}$ has exactly $n$ (not necessarily distinct) basic divisors, with each ${\cal P}^k$ counted $k$-many times. Hence to every ${\cal F}\in L$ we can adjoin a ``pattern" $\alpha_{\cal F}$: a function mapping each basic ultrafilter to its ``multiplicity" within $\cal F$. To each pattern $\alpha_{\cal F}$ corresponds a family $F_{\alpha_{\cal F}}\subseteq{\cal F}\cap{\cal U}$ determining the basic divisors of $\cal F$.

Within $L$ the $=_\sim$-equivalence classes are singletons. Let us recall (a corollary of) Theorem 3.6 from \cite{So3}, which will serve as a counterexample several times throughout this paper.

\begin{pp}\label{ramsey}
There is ${\cal P}\in\overline{\mathbb{P}}\setminus\mathbb{P}$ such that, for every $n>1$, there are at least two ultrafilters ${\cal F},{\cal G}\in\overline{L_n}$ distinct from ${\cal P}^n$ having $\cal P$ as their only basic divisor.
\end{pp}

However, in $\beta\mathbb{N}\setminus L$ things look differently: the relation $\widemid$ is not well-founded in general, and ultrafilters can not be organized in levels. The (admittedly pretentious) title of this paper addresses the possibility to nevertheless fragment (in the same sense as above) each ultrafilter into basic ones. For this purpose we will use methods of nonstandard arithmetic, introduced into the context of ultrafilter divisibility in \cite{So4}.

Probably the simplest way to introduce nonstandard extensions is via superstructures. Let $X$ be a set containing (a copy of) $\mathbb{N}$. We assume that elements of $X$ are atoms: none of them contains in its transitive closure any of the others. Let $V_0(X)=X$, $V_{n+1}(X)=V_n(X)\cup P(V_n(X))$ for $n\in\omega$ and $V(X)=\bigcup_{n<\omega}V_n(X)$. Then $V(X)$ is called a {\it superstructure}, and its {\it nonstandard extension} is a pair $(V(Y),*)$, where $V(Y)$ is a superstructure with the set of atoms $Y\supseteq X$ and $*:V(X)\rightarrow V(Y)$ is a function such that $\zve X=Y$ and satisfying\\

{\it The Transfer Principle.} For every bounded formula $\varphi(x_1,\dots,x_n)$ and every $a_1,a_2,\dots$, $a_n\in V(X)$,
$$V(X)\models\varphi(a_1,a_2,\dots,a_n)\mbox{ if and only if }V(Y)\models\varphi(\zve a_1,\zve a_2,\dots,\zve a_n).$$ 

Here, a first-order formula is bounded if all its quantifiers are bounded, i.e.\ of the form $(\forall x\in y)$ or $(\exists x\in y)$. The values of free variables $x_1,x_2,\dots,x_n$ that appear in $\varphi$ on the left can be any $a_1,a_2,\dots,a_n\in V(X)$ and on the right they are replaced with $\zve a_i\in V(Y)$. The atomic subformulas in $\varphi$ are of the form $A(x_1,\dots,x_k)$ for some $k$-ary relation $A\in V(X)$, which also gets replaced with $\zve A$. However, it is a common practice to omit $\zve{}$ in front of $\in$, arithmetic operations, the standard order $\leq$ on $\zve{\mathbb{N}}$, the operation of making $n$-tuples etc.\ to make formulas easier to read. By $\zve\varphi(x_1,x_2,\dots,x_n)$ we will denote the star-counterpart of $\varphi(x_1,x_2,\dots,x_n)$, obtained by adding stars to all values of free variables except those explicitly stated ($x_1,x_2,\dots,x_n$).

Consider a few examples that we will make use of later. First, $\zve(A^n)=(\zve A)^n$: if we denote $B=A^n$, then $(\forall a_1,\dots,a_n\in A)(a_1,\dots,a_n)\in B$ implies $(\forall a_1,\dots,a_n\in\zve A)(a_1,\dots,a_n)\in\zve B$ and $(\forall b\in B)(\exists a_1,\dots,a_n\in A)b=(a_1,\dots,a_n)$ implies $(\forall b\in\zve B)(\exists a_1,\dots,a_n\in\zve A)b=(a_1,\dots,a_n)$. In a similar way we see that $\zve(A^{(n)})=(\zve A)^{(n)}$ and $\zve(A\hspace{-0.1cm}\uparrow)=(\zve A)\hspace{-0.1cm}\uparrow$, so in all these cases we can omit the parenthesis. Also, if ${\rm dom}(f)=A\subseteq\mathbb{N}$, then ${\rm dom}(\zve f)=\zve A$. Finally, Transfer easily implies that $(\zve{\mathbb{N}}\setminus\mathbb{N},\leq)$ is divided into blocks order-isomorphic to $(\mathbb{Z},\leq)$ which Robinson, the father of modern nonstandard analysis, called galaxies in his groundbreaking book \cite{R}.

There is a well-known connection between nonstandard and ultrafilter extensions of $\mathbb{N}$: for every $x\in\zve{\mathbb{N}}$,
$$tp(x/\mathbb{N}):=\{A\subseteq\mathbb{N}:x\in\zve A\}$$
is an ultrafilter. $x$ is then called a generator of that ultrafilter. This agrees with extensions of functions $f:\mathbb{N}\rightarrow\mathbb{N}$: $tp(\zve f(x)/\mathbb{N})=\widetilde{f}(tp(x/\mathbb{N}))$ for every $x\in\zve{\mathbb{N}}$.

\begin{ex}\label{exdifprod}
Clearly, prime nonstandard numbers $p$ are exactly those belonging to $\zve{\mathbb{P}}$, so their respective ultrafilters are exactly $tp(p/\mathbb{N})={\cal P}\in\overline{\mathbb{P}}$. Also, their powers $p^k$ (for $k\in\mathbb{N}$) correspond to basic ultrafilters ${\cal P}^k$.

The best way to understand the distinction between ${\cal P}^2$ and ${\cal P}\cdot{\cal P}$ is to see that the generators of ${\cal P}\cdot{\cal P}$ are products $pq$ of distinct generators $p$ and $q$ of ${\cal P}$: Transfer easily implies that $pq\in \zve{\mathbb{P}^{(2)}}$, a set disjoint from the set $\zve{\mathbb{P}^2}$ of squares of primes.

If $tp(p/\mathbb{N})=tp(p'/\mathbb{N})={\cal P}$ and $tp(q/\mathbb{N})=tp(q'/\mathbb{N})={\cal Q}$ are distinct prime ultrafilters, $tp(pq/\mathbb{N})$ and $tp(p'q'/\mathbb{N})$ need not be the same. Namely, each of them may equal ${\cal P}\cdot{\cal Q}$ (if $(p,q)$ is a tensor pair, see below) or ${\cal Q}\cdot{\cal P}$ (if $(q,p)$ is a tensor pair), which is usually not the same, or something different from both of them.

This also explains the existence of ultrafilters $\cal P$ as described in Proposition \ref{ramsey}. A similar situation occurs at higher levels of $L$.
\end{ex}

Atoms and sets belonging to some $\zve A$ for $A\in V(X)$ are called {\it internal}. Hence, one must be careful: quantifiers of the form $(\forall x\in\zve A)$ or $(\exists x\in\zve A)$ refer only to internal sets. By {\it The Internal Definition Principle}, sets defined from internal sets are also internal; for a precise formulation see \cite{G2}, Section 13.15. For every hyperfinite set $S$ (that is, for $S\in\zve([\mathbb{N}]^{<\aleph_0})$, where $[\mathbb{N}]^{<\aleph_0}$ is the family of all finite subsets of $\mathbb{N}$) there is a unique $t\in\zve{\mathbb{N}}$ for which an internal bijection $f:\{1,2,\dots,t\}\rightarrow S$ exists; this $t$ is called the internal cardinality of $S$.

Nonstandard extensions are not unique, and they differ by richness of objects they contain. A nonstandard extension is ${\goth c}^+$-{\it saturated} if, for every family $F$ of internal sets with the finite intersection property such that $|F|\leq{\goth c}$, there is an element in $\bigcap F$. {\bf We will assume all the time that we are working with a fixed ${\goth c}^+$-saturated extension in which every initial segment of the form $|\{x\in\mathbb{N}:x\leq z\}|$ for $z\in\zve{\mathbb{N}}\setminus\mathbb{N}$ has the same cardinality.} This cardinality will be denoted by $\infty$, and we will let $\mathbb{N}_\infty=\omega\cup\{\infty\}$. This condition does not require an additional set theoretic assumption: its second half is fulfilled in nonstandard universes with property $\Delta_1$, introduced in \cite{DN2} (see Corollary 2.3 there), and the conjuction of $\Delta_1$ and ${\goth c}^+$-saturation is, by Corollary 3.6 of \cite{DN2}, equivalent to Henson's ${\goth c}^+$-isomorphism property introduced in \cite{He}.\\

With the assumption of ${\goth c}^+$-saturation, for every ${\cal F}\in\beta\mathbb{N}$ its monad
$$\mu({\cal F}):=\{x\in\zve{\mathbb{N}}:tp(x/\mathbb{N})={\cal F}\}$$
is nonempty. It also implies a connection between divisibility relations $\widemid$ and $\zvemid$, as shown in the next result (part of Theorem 3.4 from \cite{So5}).

\begin{pp}\label{ekviv}
The following conditions are equivalent for every two ultrafilters ${\cal F},{\cal G}\in\beta \mathbb{N}$:

(i) ${\cal F}\widemid{\cal G}$;

(ii) there are $x\in\mu({\cal F})$ and $y\in\mu({\cal G})$ such that $x\zvemid y$;

(iii) for every $x\in\mu({\cal F})$ there is $y\in\mu({\cal G})$ such that $x\zvemid y$;

(iv) for every $y\in\mu({\cal G})$ there is $x\in\mu({\cal F})$ such that $x\zvemid y$.
\end{pp}

Let $A\in{\cal U}$ be arbitrary, let $B$ be the set of $\mid$-minimal elements of $A$ and $B_n=B\cap L_n$. An ultrafilter that contains some $B_n$ belongs to $\overline{L_n}$; this is exactly the kind of ultrafilters studied in \cite{So3}. Note that $A=\bigcup_{n\in\omega}B_n\hspace{-1mm}\uparrow$. Let $P_A=\{p\in\mathbb{P}:(\exists n\in\mathbb{N})p^n\in A\}$ and define the function $h_A:P_A\rightarrow\mathbb{N}$ as follows:
\begin{equation}\label{eqhA}
h_A(p)=\min\{n\in\mathbb{N}:p^n\in A\}.
\end{equation}
By Transfer, for all $p\in\zve P_A$ and all $x\in\zve{\mathbb{N}}$,
\begin{equation}\label{eqhA2}
p^x\in\zve A\mbox{ if and only if }x\geq\zve h_A(p).
\end{equation}
For $p\in\zve{\mathbb{P}}\setminus\zve P_A$, no power of $p$ in is $\zve A$. Thus, in order for $tp(p^x/\mathbb{N})\widemid{\cal F}$ to hold, ${\cal F}$ needs to contain all $A\in{\cal U}$ for which $\zve h_A(p)\leq x$. In the reverse direction, for any $A\subseteq\mathbb{P}$ and any $h:A\rightarrow\mathbb{N}$ we can define
\begin{equation}\label{eqhA3}
A^h=\{m\in\mathbb{N}:(\exists p\in A)p^{h(p)}\mid m\}\in{\cal U}.
\end{equation}

If $\{{\cal F}_i:i\in I\}$ is a set of ultrafilters and ${\cal W}$ an ultrafilter on $I$, ${\cal G}=\lim_{i\rightarrow{\cal W}}{\cal F}_i$ is the ultrafilter defined by: $A\in{\cal G}$ if and only if $\{i\in I:A\in{\cal F}_i\}\in{\cal W}$. The following was proven in \cite{So5}, Lemma 4.1.

\begin{pp}\label{sublimit}
(a) Every chain $\langle[{\cal F}_i]_\sim:i\in I\rangle$ in $(\beta \mathbb{N}/\hspace{-1mm}=_\sim,\widemid)$ has the least upper bound $[{\cal G}_U]_\sim$ (obtained as ${\cal G}_U=\lim_{i\rightarrow{\cal W}}{\cal F}_i$ for any ${\cal W}$ containing all final segments of $I$) and the greatest lower bound $[{\cal G}_L]_\sim$.

(b) $\bigcup_{i\in I}({\cal F}_i\cap{\cal U})={\cal G}_U\cap{\cal U}$ and $\bigcap_{i\in I}({\cal F}_i\cap{\cal U})={\cal G}_L\cap{\cal U}$.
\end{pp}

Since the $=_\sim$-class of the l.u.b.\ ${\cal G}_U$ does not depend on the choice of $\cal W$, in the case when $I=\gamma$ is an ordinal we will denote by $\lim_{\delta\rightarrow\gamma}{\cal F}_\delta$ the class $[{\cal G}_U]_\sim$.

In any nonstandard extension the following generalization of the Fundamental theorem of arithmetic holds. It is easily obtained from the Transfer principle. Let $\langle p_n:n\in\mathbb{N}\rangle$ be the increasing enumeration of $\mathbb{P}$, so its nonstandard extension $\langle p_n:n\in\zve{\mathbb{N}}\rangle$ is the increasing enumeration of $\zve{\mathbb{P}}$. Recall that, for $p\in\mathbb{P}$, $n,k\in\mathbb{N}$, $p^k\parallel n$ means that $k$ is the largest natural number $l$ such that $p^l\mid n$; we say that $p^k$ is an exact divisor of $n$. Likewise, for $p\in\zve{\mathbb{P}}$, $x,k\in\zve{\mathbb{N}}$, $p^k\zvepar x$ means that $k$ is the largest $l\in\zve{\mathbb{N}}$ such that $p^l\zvemid x$. If $p^k\zvepar x$, we also write $k={\rm exp}_px$.

\begin{pp}\label{fund}
(a) For every $z\in\zve{\mathbb{N}}$ and every internal sequence $\langle h(n):n\leq z\rangle$ there is a unique $x\in\zve{\mathbb{N}}$ such that $p_n^{h(n)}\zvepar x$ for $n\leq z$ and $p_n\nzvemid x$ for $n>z$; we denote this element by $\prod_{n\leq z}p_n^{h(n)}$.

(b) Every $x\in\zve{\mathbb{N}}$ can be uniquely represented as $\prod_{n\leq z}p_n^{h(n)}$ for some $z\in\zve{\mathbb{N}}$ and some internal sequence $\langle h(n):n\leq z\rangle$ such that $h(z)>0$.
\end{pp}

The product ${\cal F}\cdot{\cal G}$ is a $\widemid$-minimal ultrafilter divisible by both $\cal F$ and $\cal G$. However, as mentioned in Example \ref{exdifprod}, it may be only one of many such ultrafilters. More precisely, by Lemma 3.7 from \cite{So3}, they are exactly ulrafilters containing the filter generated by sets $AB$ for $A\in{\cal F}$, $B\in{\cal G}$. Thus, products $xy$ for $x\in\mu({\cal F})$ and $y\in\mu({\cal G})$ do not always belong to ${\cal F}\cdot{\cal G}$. They do whenever $(x,y)$ is a {\it tensor pair}, meaning that it belongs to $\mu({\cal F}\otimes{\cal G})$, where ${\cal F}\otimes{\cal G}=\{S\subseteq\mathbb{N}:\{x\in\mathbb{N}:\{y\in\mathbb{N}:(x,y)\in S\}\in{\cal G}\}\in{\cal F}\}$ is the tensor product of ultrafilters $\cal F$ and $\cal G$. Theorem 3.4 from \cite{P} gives an equivalent condition that we will use (and several more can be found in Theorem 11.5.7 of \cite{DN}).

\begin{pp}\label{tensor}
$(x,y)$ is a tensor pair if and only if, for every $f:\mathbb{N}\rightarrow\mathbb{N}$, either $\zve f(y)\in\mathbb{N}$ or $\zve f(y)>x$.
\end{pp}

The following lemma is a version of Theorem 11.5.12 from \cite{DN}.
 
\begin{lm}\label{tensor2}
Let ${\cal F},{\cal G}\in\beta\mathbb{N}$ and ${\cal H}\in\beta(\mathbb{N}\times\mathbb{N})$. If $x_0\in\mu({\cal F})$ and $y_0\in\mu({\cal G})$ are such that $(x_0,y_0)\in\mu({\cal H})$, then for every $x_1\in\mu({\cal F})$ there is $y_1\in\mu({\cal G})$ such that $(x_1,y_1)\in\mu({\cal H})$.

In particular, for every $x\in\zve{\mathbb{N}}\setminus\mathbb{N}$ and every ${\cal G}\in\beta\mathbb{N}\setminus\mathbb{N}$, there is $y\in\mu({\cal G})$ such that $(x,y)$ is a tensor pair; analogously there is $y'\in\mu({\cal G})$ such that $(y',x)$ is a tensor pair.
\end{lm}

\dokaz Let $B_{A,X}=\{y\in\zve X:(x_1,y)\in\zve A\}$ for $A\subseteq\mathbb{N}\times\mathbb{N}$ and $X\subseteq\mathbb{N}$. Let $F=\{B_{A,X}:A\in{\cal H}\land X\in{\cal G}\}$. Since $F$ is closed for finite intersections, to prove that $F$ has the finite intersection property it suffices to show that $B_{A,X}\neq\emptyset$ for all $A\in{\cal H}$, $X\in{\cal G}$. If we denote by $\pi_1$ the first projection, $y_0$ witnesses that $x_0\in\zve\pi_1[\zve A\cap(\zve{\mathbb{N}}\times\zve X)]=\zve(\pi_1[A\cap(\mathbb{N}\times X)])$, so $x_1\in\zve\pi_1[\zve A\cap(\zve{\mathbb{N}}\times\zve X)]$ as well. This means that there is $y\in\zve X$ such that $(x_1,y)\in\zve A$. Now ${\goth c}^+$-saturation implies that there is $y_1\in\bigcap F$, which means that $y_1\in\mu({\cal G})$ and $(x_1,y_1)\in\mu({\cal H})$.\kraj

\section{Basic ultrafilters}

We begin the description of the divisibility order by describing powers of prime ultrafilters. If $p\in\zve{\mathbb{P}}$ is a fixed nonstandard prime and ${\cal P}=tp(p/\mathbb{N})$, $\langle p^x:x\in\zve{\mathbb{N}}\rangle$ is a $\zvemid$-increasing chain in $\zve{\mathbb{N}}$, and so $\langle tp(p^x/\mathbb{N}):x\in\zve{\mathbb{N}}\rangle$ is a $\widemid$-nondecreasing chain in $\beta\mathbb{N}$. For distinct $m,n\in\mathbb{N}$ the ultrafilters ${\cal P}^m=tp(p^m/\mathbb{N})$ and ${\cal P}^n=tp(p^n/\mathbb{N})$ are $=_\sim$-nonequivalent. What about exponents from $\zve{\mathbb{N}}\setminus\mathbb{N}$? Example 4.5 from \cite{So4} shows that for $p\in\mathbb{P}$ the situation is simple.

\begin{pp}\label{caseP}
If $p\in\mathbb{P}$, all the ultrafilters $tp(p^x/\mathbb{N})$ for $x\in\zve{\mathbb{N}}\setminus\mathbb{N}$ are $=_\sim$-equivalent.
\end{pp}

However, for $p\in\zve{\mathbb{P}}\setminus\mathbb{P}$ this need not be true, as the next example shows.

\begin{ex}\label{cut}
Again, let $\langle p_x:x\in\zve{\mathbb{N}}\rangle$ be the increasing enumeration of $\zve{\mathbb{P}}$, and let $A=\bigcup_{n\in\mathbb{N}}\{{p_n}^n\}\hspace{-0.1cm}\uparrow$. (In terms of (\ref{eqhA3}): if $h:\mathbb{P}\rightarrow\mathbb{N}$ is given by $h(p_n)=n$, then $A=A^h$.) Then:

(a) For every $x\in\zve{\mathbb{N}}$, Transfer implies that $p_x^x\in\zve A$ but $p_x^{x-1}\notin\zve A$. Since $A\in{\cal U}$, this means that the consecutive powers $p_x^{x-1}$ and $p_x^x$ are in the monads of $=_\sim$-nonequivalent ultrafilters.


(b) $(\forall p,q\in\mathbb{P})(p\neq q\Rightarrow(\exists n\in\mathbb{N})(p^n\in A\;\underline{\lor}\; q^n\in A))$ is also true (where $\underline{\lor}$ is the exclusive disjunction), and therefore
$$(\forall p,q\in\zve{\mathbb{P}})(p\neq q\Rightarrow(\exists x\in\zve{\mathbb{N}})(p^x\in\zve A\;\underline{\lor}\; q^x\in\zve A)).$$
Thus, any two prime nonstandard numbers have powers $p^x$ and $q^x$ such that $tp(p^x/\mathbb{N})\neq_\sim tp(q^x/\mathbb{N})$.
\end{ex}

\begin{de}
Let ${\cal P}\in\overline{\mathbb{P}}$. The relation $\approx_{\cal P}$ on $\mu({\cal P})\times\zve{\mathbb{N}}$ is defined as follows:
$$(p,x)\approx_{\cal P}(q,y)\mbox{ if and only if }tp(p^x/\mathbb{N})=_\sim tp(q^y/\mathbb{N}).$$
$\approx_{\cal P}$ is, of course, an equivalence relation; let
$${\cal E}_{\cal P}=\{[(p,x)]_{\approx_{\cal P}}:(p,x)\in\mu({\cal P})\times\zve{\mathbb{N}}\}$$
be the set of its equivalence classes. For any $p\in\mu({\cal P})$ and $u\in{\cal E}_{\cal P}$, let the ``vertical section" be the set defined by $u_p:=\{x:(p,x)\in u\}$.
\end{de}

Example \ref{cut}(b) shows that $tp(p^x/\mathbb{N})=_\sim tp(q^x/\mathbb{N})$ need not hold for $p,q\in\mu({\cal P})$, so the sets $u_p$ are not independent of the choice of $p$. 

\begin{de}
Families of ultrafilters of the form ${\cal P}^u:=\{tp(p^x/\mathbb{N}):(p,x)\in u\}$ for some ${\cal P}\in\overline{\mathbb{P}}$ and $u\in{\cal E}_{\cal P}$ will be called basic. Also, let $\mu({\cal P}^u)=\bigcup_{{\cal F}\in{\cal P}^u}\mu({\cal F})$.

By ${\cal B}$ we denote the set of all basic classes.
\end{de}

We will think of the $=_\sim$-classes ${\cal P}^u$ as powers of ${\cal P}$.

\begin{lm}\label{approxosob}
Let ${\cal P}\in\overline{\mathbb{P}}$ and $u\in{\cal E}_{\cal P}$.

(a) All elements of $\mu({\cal P}^u)$ are of the form $p^x$ for some $(p,x)\in u$.

(b) For every $p\in\mu({\cal P})$ the set $u_p$ is nonempty and convex: if $x,y\in u_p$ and $x<z<y$, then $z\in u_p$ as well.

(c) Each $u_p$ is either a singleton or a union of galaxies.
\end{lm}

\dokaz (a) If ${\cal F}\in{\cal P}^u$, there is $(q,y)\in u$ such that $q^y\in\mu({\cal F})$. $q^y\in\zve{\mathbb{P}}^{exp}$, so $\mathbb{P}^{exp}\in{\cal F}$. Thus, every element of $\mu({\cal F})$ belongs to $\zve{\mathbb{P}}^{exp}$, and so it is of the form $p^x$. Clearly $(p,x)\approx_{\cal P}(q,y)$, so $(p,x)\in u$.

(b) To prove $u_p\neq\emptyset$ let $(q,y)\in u$ and ${\cal F}=tp(q^y/\mathbb{N})$. Using Proposition \ref{ekviv}, $q\zvemid q^y$ implies ${\cal P}\widemid{\cal F}$, which in turn implies that there is an element of $\mu({\cal F})$ divisible by $p$, which by (a) must be of the form $p^x$ for some $x\in\zve{\mathbb{N}}$. Hence $(p,x)\in u$. The convexity of $u_p$ is obvious.

(c) Assume that, for some $x\in\zve{\mathbb{N}}\setminus\mathbb{N}$, $tp(p^x/\mathbb{N})\neq_\sim tp(p^{x+1}/\mathbb{N})$; let us prove that $tp(p^y/\mathbb{N})\neq_\sim tp(p^{y+1}/\mathbb{N})$ holds for every element $y$ of the galaxy of $x$. By the assumption there is $A\in{\cal U}\cap tp(p^{x+1}/\mathbb{N})\setminus tp(p^x/\mathbb{N})$, in other words $\zve h_A(p)=x+1$ (see (\ref{eqhA2})). Let $g:\mathbb{N}\setminus\{1\}\rightarrow\mathbb{N}$ be the function defined by $g(p_1^{a_1}\dots p_k^{a_k})=p_1^{a_1+1}\dots p_k^{a_k+1}$ (for distinct $p_1,\dots,p_k\in\mathbb{P}$ and $a_i>0$). Then $B:=g[A]\hspace{-0.1cm}\uparrow\in{\cal U}$. $\zve g(p^{x+1})=p^{x+2}\in\zve B$ and $p^{x+1}\notin\zve B$, so $\zve h_{B}(p)=x+2$. This implies that $tp(p^{x+1}/\mathbb{N})\neq_\sim tp(p^{x+2}/\mathbb{N})$ and inductivelly $tp(p^{x+n}/\mathbb{N})\neq_\sim tp(p^{x+n+1}/\mathbb{N})$. Using $f:\mathbb{N}\setminus\{1\}\rightarrow\mathbb{N}$ defined by $f(p_1^{a_1}\dots p_k^{a_k})=p_1^{a_1-1}\dots p_k^{a_k-1}$, in a similar way we conclude $tp(p^{x-n}/\mathbb{N})\neq_\sim tp(p^{x-n+1}/\mathbb{N})$ for $n\in\mathbb{N}$.\kraj

\begin{de}
On ${\cal E}_{\cal P}$ we define the relation:
\begin{eqnarray*}
u\prec_{\cal P}v &\mbox{ if and only if } & u\neq v\mbox{ and for some } p\in\mu({\cal P})\mbox{ and some }x,y\in\zve{\mathbb{N}}\\
&& \mbox{ holds }(p,x)\in u,(p,y)\in v\mbox{ and }x<y.
\end{eqnarray*}
We write $u\preceq_{\cal P}v$ if $u\prec_{\cal P}v$ or $u=v$.
\end{de}

Proposition \ref{ekviv} implies that, if $u\prec_{\cal P}v$, then in fact for all $p\in\mu({\cal P})$ and all $x,y\in\zve{\mathbb{N}}$ such that $(p,x)\in u$ and $(p,y)\in v$, $x<y$ holds.

\begin{lm}\label{supinf}
For every ${\cal P}\in\overline{\mathbb{P}}$:

(a) $\prec_{\cal P}$ is a strict linear order.

(b) Every increasing sequence in $({\cal E}_{\cal P},\prec_{\cal P})$ has a supremum and every decreasing sequence has an infimum.
\end{lm}

\dokaz (a) is obvious, using Lemma \ref{approxosob}(b) and the remark preceding this lemma.

(b) We can assume, without loss of generality, that the given sequence is well-ordered (otherwise we can first thin it out into a cofinal well-ordered subsequence). So let $\langle u_\xi:\xi<\gamma\rangle$ be a $\prec_{\cal P}$-increasing sequence in ${\cal E}_{\cal P}$ of some limit length $\gamma$. Pick some $p\in\mu({\cal P})$ and $x_\xi$ for $\xi<\gamma$ such that $(p,x_\xi)\in u_\xi$. Then $\langle tp(p^{x_\xi}/\mathbb{N}):\xi<\gamma\rangle$ is a $\widemid$-increasing sequence of ultrafilters so by Proposition \ref{sublimit} it has a least upper bound ${\cal G}$. Since each $tp(p^{x_{\xi+1}}/\mathbb{N})\setminus tp(p^{x_\xi}/\mathbb{N})$ contains a set in $\cal U$ and $|\cal U|={\goth c}$, the sequence is of length less than ${\goth c}^+$. But ${\goth c}^+$-saturation easily proves that the cofinality of $\zve{\mathbb{N}}$ is at least ${\goth c}^+$, so there is an upper bound $z$ for $\{x_{\xi}:\xi<\gamma\}$. Using Proposition \ref{ekviv} we find $y\leq z$ such that $p^y\in\mu({\cal G})$, so $[(p,y)]_{\approx_{\cal P}}$ must be the supremum of the given sequence.

Analogously we prove that $({\cal E}_{\cal P},\prec_{\cal P})$ is closed for infimums of decreasing sequences.\kraj

Inspecting the proof of Theorem 4.6 from \cite{So5}, we can conclude that, for every ${\cal P}\in\overline{\mathbb{P}}\setminus\mathbb{P}$, the order $({\cal E}_{\cal P},\prec_{\cal P})$ contains a copy of $(\mathbb{R},<)$. We will not use this fact here.

Note that both cases from Lemma \ref{approxosob}(c) are possible: the first case occurs by Example \ref{cut}(a), and above every such galaxy (``cut" into singleton classes) comes, by Lemma \ref{supinf}(b), the second case: a union of galaxies. For example, if $u_p^x=\{x\}$ are singleton classes (for $x$ from some galaxy $X$), then $\sup_{x\in X}u_p^x$ can not be a singleton $\{z\}$ because $z-1$ would be a smaller upper bound.

As another corollary of Lemma \ref{supinf}, we conclude that each ${\cal E}_{\cal P}$ has the greatest element $u_{max}$; let us write ${\cal P}^{max}$ for ${\cal P}^{u_{max}}$. Also, for ${\cal P}\in\overline{\mathbb{P}}$ we will denote ${\cal P}^\omega={\cal P}^u$, where $u=\sup\omega$ in ${\cal E}_{\cal P}$. As we already mentioned, for $p\in\mathbb{P}$ holds $p^\omega=p^{max}$, so the order-type of $({\cal E}_p,\prec_p)$ is $\omega+1$. When convenient, we will identify the first $\omega+1$ many elements of ${\cal E}_{\cal P}$ with the set of ordinals $\omega+1$.

\begin{ex}
Let $A\in{\cal U}$, ${\cal P}\in\overline{\mathbb{P}}\setminus\mathbb{P}$, $x\in\zve{\mathbb{N}}\setminus\mathbb{N}$ and $p\in\mu({\cal P})$.

(a) If $(p,x)$ is a tensor pair, then $p^x$ belongs to $\mu({\cal P}^{max})$. To see that, assume the opposite: there is some $y>x$ such that $[(p,x)]_{\approx_{\cal P}}\prec_{\cal P}[(p,y)]_{\approx_{\cal P}}$. By Lemma \ref{tensor2} there is $z$ such that $(y,z)$ is a tensor pair, in particular $p<x<y<z$. But then by Proposition \ref{tensor} $(p,z)$ is also a tensor pair, so $[(p,x)]_{\approx_{\cal P}}\prec_{\cal P}[(p,y)]_{\approx_{\cal P}}\preceq_{\cal P}[(p,z)]_{\approx_{\cal P}}=[(p,x)]_{\approx_{\cal P}}$, a contradiction.

(b) If $(x,p)$ is a tensor pair, then $p^x$ belongs to $\mu({\cal P}^\omega)$. Namely, by Proposition \ref{tensor} for every $A\in{\cal U}$ either $\zve h_A(p)\in\mathbb{N}$ or $\zve h_A(p)>x$. By (\ref{eqhA2}), $p^x\in\zve A$ can hold only if $\zve h_A(p)=m$ for some $m\in\mathbb{N}$, which means that already $p^m\in\zve A$, so $A\in{\cal P}^m$. Thus, $tp(p^x/\mathbb{N})$ contains only sets from ${\cal U}$ belonging to some ${\cal P}^m$, so it is a least upper bound of $\{{\cal P}^m:m\in\mathbb{N}\}$.
\end{ex}

If $A\subseteq\mathbb{N}$ and ${\cal P}^u\in{\cal B}$, let us abuse the notation and write $A\in{\cal P}^u$ (or ${\cal P}^u\in\overline{A}$) if $A\in{\cal F}$ for all ${\cal F}\in{\cal P}^u$. For example, $\mathbb{P}^{exp}\in{\cal P}^u$ for all ${\cal P}^u\in{\cal B}$. Note that, when $A\in{\cal C}$, for $A\in{\cal P}^u$ it suffices that  $A\in{\cal F}$ for some ${\cal F}\in{\cal P}^u$. Hence we denote ${\cal P}^u\cap{\cal U}=\{A\in{\cal U}:A\in{\cal P}^u\}$ and ${\cal P}^u\cap{\cal C}=\{A\in{\cal C}:A\in{\cal P}^u\}$.

Let us call the topology on $\cal B$ generated by base sets $\overline{B}$ for $B\in{\cal U}$ the $\cal U$-topology. The closure of $S\subseteq{\cal B}$ in this topology will be denoted ${\rm cl}_{\cal U}S$.

\begin{lm}\label{utop}
Open sets of the $\cal U$-topology on $\cal B$ are exactly the sets of the form $\overline{A^h}$ for $h:A\rightarrow\mathbb{N}$ such that $A\subseteq\mathbb{P}$.
\end{lm}

\dokaz Let $B\in\cal U$; we want to find a set of the form $A^h$ such that, for ${\cal Q}^u\in{\cal B}$, $B\in{\cal Q}^u$ if and only if $A^h\in{\cal Q}^u$. Simply define $h=h_B$ as in (\ref{eqhA}). Clearly $A^{h_B}\subseteq B$, which gives us one of the desired implications. On the other hand, for any ${\cal Q}^u\in\overline{B}$, $B\cap\mathbb{P}^{exp}\in{\cal Q}^u$ so $A^{h_B}=(B\cap\mathbb{P}^{exp})\hspace{-0.1cm}\uparrow\in{\cal Q}^u$ as well.\kraj

\begin{ex}\label{zatvex0}
To every $A\subseteq\mathbb{P}$ corresponds a set $A\hspace{-0.1cm}\uparrow\in{\cal U}$ so that, for every ${\cal Q}\in\overline{\mathbb{P}}$, ${\cal Q}\in\overline{A}$ if and only if ${\cal Q}\in\overline{A\hspace{-0.1cm}\uparrow}$. Thus the $\cal U$-topology on $\overline{\mathbb{P}}$ coincides with the standard one (generated by all $\overline{A}$ for $A\subseteq\mathbb{P}$). If $S\subseteq\overline{\mathbb{P}}$ then ${\cal Q}\in{\rm cl}_{\cal U}S$ implies that, for every $k\in\mathbb{N}$, ${\cal Q}^k\in{\rm cl}_{\cal U}\{{\cal P}^k:{\cal P}\in S\}$ as well.

In fact, if ${\cal G}={\cal Q}^k\in\overline{\mathbb{P}^k}$, then essentially all sets $S'$ such that ${\cal G}\in{\rm cl}_{\cal U}S'$ are obtained in this way: for any $A\in{\cal G}\cap{\cal U}$, $A\cap(\mathbb{P}^k)\in{\cal G}$ as well, so there is an ultrafilter ${\cal P}^k\in S'\cap\overline{\mathbb{P}^k}$ such that $A\in{\cal P}^k$, meaning that ${\cal G}\in{\rm cl}_{\cal U}(S'\cap\overline{\mathbb{P}^k})$. Hence in such occasions we can immediately restrict ourselves to the case $S'=\{{\cal P}^k:{\cal P}\in S\}$ for some $S\subseteq\overline{\mathbb{P}}$, whence ${\cal G}={\cal Q}^k\in{\rm cl}_{\cal U}S'$ if and only if ${\cal Q}\in{\rm cl}_{\cal U}S$.
\end{ex}

The following lemma will provide a better understanding of the situation when a set $u_p$ is a singleton.

\begin{lm}\label{selfdivpom}
Let ${\cal P}\in\overline{\mathbb{P}}$ and $u\in{\cal E}_{\cal P}$.

(a) For any 
$f:\mathbb{P}\rightarrow\mathbb{N}$ and $p,q\in\mu({\cal P})$, $p^{\zve f(p)}\in\mu({\cal P}^u)$ if and only if $q^{\zve f(q)}\in\mu({\cal P}^u)$.

(b) If a vertical section $u_p$ contains an element $p^{\zve f(p)}$ for some $f:\mathbb{P}\rightarrow\mathbb{N}$, then it is the only element and in fact $u=\{(p,\zve f(p)):p\in\mu({\cal P})\}$.
\end{lm}

\dokaz (a) By Lemma \ref{utop} it suffices to show that, for every $A\in{\cal P}$ and every $h:A\rightarrow\mathbb{N}$, $p^{\zve f(p)}\in\zve{A^h}$ if and only if $q^{\zve f(q)}\in\zve{A^h}$. By (\ref{eqhA2}), this is equivalent to showing that $\zve f(p)\geq\zve h(p)$ if and only if $\zve f(q)\geq\zve h(q)$. If we define $X=\{b\in A:f(b)\geq h(b)\}$, then $p\in\zve X$ if and only if $q\in\zve X$, which proves the claim.

(b) From (a) it follows that all the pairs $(p,\zve f(p))$ are in $u$. But, for any such pair, $(p,\zve f(p)-1)\notin u$ because $p^{\zve f(p)-1}\notin\mathbb{P}^f$. Likewise, $(p,\zve f(p)+1)\notin u$ because $p^{\zve f(p)+1}\in\mathbb{P}^g$ (where $g(p)=f(p)+1$ for all $p\in\mathbb{P}$), while $p^{\zve f(p)}\notin\mathbb{P}^g$.\kraj

\section{Patterns}

As we saw in the previous section, when generalizing from $L$ to the whole of $\beta\mathbb{N}$ the set of basic divisors needs to be expanded by classes ${\cal P}^u$ for $u\in{\cal E}_{\cal P}\setminus\omega$. Another difference is that an ultrafilter can be divisible by a basic class infinitely many times. It will turn out that, under the assumptions stated in the introduction, there is only one possibility for such infinite multiplicity of a basic divisor: the cardinality $\infty$ of all sets of the form $\{1,2,\dots,z\}$ for $z\in\zve{\mathbb{N}}\setminus\mathbb{N}$.

\begin{de}
For $x\in\zve{\mathbb{N}}$ and $u,v\in{\cal E}_{\cal P}$ such that $u\preceq_{\cal P}v$, denote
$$D^{[u,v]}_x:=\{(p,k):(u\preceq_{\cal P}[(p,k)]_{\approx_{\cal P}}\preceq_{\cal P}v)\land p^{k}\zvepar x\}.$$
\end{de}

\begin{lm}\label{imacplus}
Let ${\cal P}\in\overline{\mathbb{P}}$, $u,v\in{\cal E}_{\cal P}$ and ${\cal F}\in\beta\mathbb{N}$. If there is $x_0\in\mu({\cal F})$ such that $D^{[u,v]}_{x_0}$ is infinite then, for every $z\in\zve{\mathbb{N}}\setminus\mathbb{N}$, there is $x\in\mu({\cal F})$ such that $D^{[u,v]}_x$ has a hyperfinite subset of internal cardinality $z$, and thus $|D^{[u,v]}_x|=\infty$.
\end{lm}

\dokaz The formula
$$\theta_1(y,C)\equiv(\forall p\in\mathbb{P})(p\mid y\Rightarrow p^{{\rm exp}_py}\in C)$$
claims that all powers of primes which are exact factors of $y$ belong to $C$.
$$\theta_2(y,t)\equiv(\exists f:\{1,2,\dots,t\}\rightarrow\mathbb{P})(f\mbox{ is one-to-one }\land (\forall i\leq t)f(i)\mid y)$$
claims that $y$ has ``at least $t$" prime divisors. Also let
$$\theta_3(x,y)\equiv (\forall p\in\mathbb{P})(p\mid y\Rightarrow {\rm exp}_px={\rm exp}_py).$$
Choose any $z\in\zve{\mathbb{N}}\setminus\mathbb{N}$ and let $B_{A,C}=\{(x,y)\in\zve A\times\zve{\mathbb{N}}:\zve\theta_1(y,\zve C)\land\zve\theta_2(y,z)\land\zve\theta_3(x,y)\}$. $B_{A,C}$ is internal, since it is defined from $\zve A$, $\zve C$ and $z$, all of which are internal. Finally, let $F=\{B_{A,C}:A\in{\cal F}\cap{\cal C}\land C\in{\cal P}^u\cap{\cal P}^v\cap{\cal C}\}$. (Recall that $\cal C$ is the family of all convex sets.)

Let $A\in{\cal F}\cap{\cal C}$ and $C\in{\cal P}^u\cap{\cal P}^v\cap{\cal C}$. Consider the formula
$$\psi\equiv(\forall t\in\mathbb{N})(\exists x\in A)(\exists y\in\mathbb{N})(\theta_1(y,C)\land\theta_2(y,t)\land\theta_3(x,y)).$$
For any $t\in\mathbb{N}$ and the formula
$$\varphi\equiv(\exists x\in A)(\exists q_1,\dots,q_t\in\mathbb{P}^{exp})((\forall i\neq j)q_i\neq q_j\land(\forall i)(q_i\parallel x\land q_i\in C))$$
its star-counterpart $\zve\varphi$ holds, as witnessed by $x_0$ and any $t$ of its exact divisors $p_i^{k_i}$ such that $u\preceq_{\cal P}[(p_i,k_i)]_{\approx_{\cal P}}\preceq_{\cal P}v$. Therefore $\varphi$ is also true, and so is $\psi$. $\zve\psi$ (for $t=z$) establishes that $B_{A,C}\neq\emptyset$. Since $F$ is closed for finite intersections, it has the finite intersection property, so by ${\goth c}^+$-saturation there are $x\in\mu({\cal F})$ and an internal one-to-one function $f:\{1,2,\dots,z\}\rightarrow\zve{\mathbb{P}}$ so that, for every $i\leq z$, $(f(i),{\rm exp}_{f(i)}x)\in D^{[u,v]}_x$. Finally, the set $\{f(i):i\leq z\}$ has internal cardinality $z$, so $|D^{[u,v]}_x|\geq\infty$. However, by Proposition \ref{fund}, the cardinality of $D^{[u,v]}_x$ can not be larger than $\infty$, so $|D^{[u,v]}_x|=\infty$.\kraj

\begin{te}\label{samocplus}
If $x\in\zve{\mathbb{N}}$, ${\cal P}\in\overline{\mathbb{P}}$ and $u,v\in{\cal E}_{\cal P}$, then $x$ has either finitely many or $\infty$-many exact divisors from $\bigcup_{u\preceq_{\cal P}w\preceq_{\cal P}v}\mu({\cal P}^w)$.
\end{te}

\dokaz Assume that $x$ has infinitely many divisors from $\bigcup_{u\preceq_{\cal P}w\preceq_{\cal P}v}\mu({\cal P}^w)$ and let ${\cal F}=tp(x/\mathbb{N})$. Choose any $z\in\zve{\mathbb{N}}\setminus\mathbb{N}$. By Lemma \ref{imacplus} there are $y\in\mu({\cal F})$ and internal one-to-one function $f:\{1,2,\dots,z\}\rightarrow D^{[u,v]}_y$. Lemma \ref{tensor2} provides some $t\in\zve{\mathbb{N}}\setminus\mathbb{N}$ such that $tp(y,z/\mathbb{N})=tp(x,t/\mathbb{N})$. For $A\subseteq\mathbb{N}$ let
$$B_A:=\{(m,k)\in\mathbb{N}^2:(\exists f:\{1,2,\dots,k\}\rightarrow A\cap\mathbb{P}^{exp})(f\mbox{ is one-to-one}\land(\forall i\leq k)f(i)\parallel m)\}.$$
Now $(y,z)\in\zve B_A$ implies $(x,t)\in\zve B_A$ for all $A\in{\cal P}^u\cap{\cal P}^v\cap{\cal C}$. Hence the family
$$\{\{f:\{1,2,\dots,t\}\rightarrow\zve(A\cap\mathbb{P}^{exp}) | f\mbox{ is one-to-one}\land(\forall i\leq t)f(i)\zvepar x\}:A\in{\cal P}^u\cap{\cal P}^v\cap{\cal C}\}$$
has the f.i.p. If $f:\{1,2,\dots,t\}\rightarrow\zve{\mathbb{P}}^{exp}$ is in its intersection, then $f(i)$ for $i\leq t$ are $\infty$-many divisors of $x$ from $\bigcup_{u\preceq_{\cal P}w\preceq_{\cal P}v}\mu({\cal P}^w)$.\kraj

There are two particularly important corollaries of the result above. We get the first by putting $v=max$: the number of exact divisors from $\bigcup_{w\succeq_{\cal P}u}\mu({\cal P}^w)$ is either finite or $\infty$ for every $u\in{\cal E}_{\cal P}$. The second one is obtained for $u=v$: every $x$ has either finitely many or $\infty$-many exact divisors from $\mu({\cal P}^u)$. This will simplify considerably the following definition of pattern of $x$.

\begin{de}\label{defpattern}
Let ${\cal A}$ be the set of all functions $\alpha:{\cal B}\rightarrow\mathbb{N}_\infty$ such that $\sum_{k\in\omega+1}\alpha(p^k)\leq 1$ for $p\in\mathbb{P}$. Elements $\alpha\in{\cal A}$ are called patterns.

If ${\cal P}^u=\{{\cal P}^k\}$ is a singleton, we identify ${\cal P}^u$ with ${\cal P}^k$ and write $\alpha({\cal P}^k)$ instead of $\alpha({\cal P}^u)$. In particular, for $k\in\mathbb{N}$, $\overline{\mathbb{P}^k}$ is regarded as a subset of $\cal B$.
\end{de}

These patterns are tailored to represent the quantities of each of basic divisors of a given ultrafilter, which will become clear in Definitions \ref{defalpha} and \ref{alphault}. The additional condition for $p\in\mathbb{P}$ will be explained in Example \ref{exalphault}(a).

We will also write $\alpha=\{({\cal P}_i^{u_i},n_i):i\in I\}$ (meaning that $\alpha({\cal P}_i^{u_i})=n_i$), omitting some (or all) of the pairs $({\cal P}^u,m)$ for which $m=0$.

\begin{de}
We will say that $\alpha\in{\cal A}$ is $\cal U$-closed if, whenever ${\cal Q}^u\in{\cal B}$ and $n\in\mathbb{N}$, if for every $A\in{\cal Q}^u\cap{\cal U}$ holds $\sum_{{\cal P}^v\in\overline{A}}\alpha({\cal P}^v)\geq n$, then $\sum_{w\succeq_{\cal Q}u}\alpha({\cal Q}^w)\geq n$.

The family of all $\cal U$-closed patterns will be denoted by ${\cal A}_{cl}$.
\end{de}

Intuitively, ${\cal U}$-closedness of a pattern $\alpha$ means that, if $\sum_{w\succeq_{\cal Q}u}\alpha({\cal Q}^w)$ is finite, then there is a neighborhood $\overline{A}$ of ${\cal Q}^u$ in which there are no basic classes ``appearing" in $\alpha$ other than ${\cal Q}^u$ or higher powers of $\cal Q$.

Some special cases should help to illuminate the concept of $\cal U$-closedness.

\begin{ex}\label{zatvex}
Let $\alpha\in{\cal A}_{cl}$.

(a) Let $S\subseteq\overline{\mathbb{P}}$, ${\cal Q}\in{\rm cl}_{\cal U}S$ and $k\in\mathbb{N}$. Recall from Example \ref{zatvex0} that ${\cal Q}^k$ is in the $\cal U$-closure of $\{{\cal P}^k:{\cal P}\in S\}$. If, for every such $S$, $\sum_{{\cal P}\in S}\sum_{u\succeq_{\cal P}k}\alpha({\cal P}^u)\geq n$, then $\sum_{u\succeq_{\cal Q}k}\alpha({\cal Q}^u)\geq n$. In particular, $\sum_{{\cal P}\in S}\sum_{u\succeq_{\cal P}k}\alpha({\cal P}^u)=\infty$ for every such $S$ implies that $\sum_{u\succeq_{\cal Q}k}\alpha({\cal Q}^u)=\infty$.

(b) Let $n\in\mathbb{N}$. If $\langle v_\xi:\xi<\gamma\rangle$ is a $\prec_{\cal P}$-increasing sequence in ${\cal E}_{\cal P}$ and $u=\sup_{\xi<\gamma}v_\xi$ (see Lemma \ref{supinf}(b)), then ${\cal P}^u\in{\rm cl}_{\cal U}(\{{\cal P}^{v_\xi}:\xi<\gamma\})$, so $\sum_{w\succeq_{\cal P}v_\xi}\alpha({\cal P}^w)\geq n$ for all $\xi<\gamma$ implies $\sum_{w\succeq_{\cal P}u}\alpha({\cal P}^w)\geq n$. Consequently, $\sum_{w\succeq_{\cal P}v_\xi}\alpha({\cal P}^w)=\infty$ for $\xi<\gamma$ implies $\sum_{w\succeq_{\cal P}u}\alpha({\cal P}^w)=\infty$.
\end{ex}

To every $\alpha\in{\cal A}$ and every prime ${\cal P}$ we can adjoin a sequence $\alpha\upharpoonright{\cal P}:=\langle\alpha({\cal P}^u):u\in{\cal E}_{\cal P}\rangle$. Clearly, fixing $\alpha\upharpoonright{\cal P}$ for every ${\cal P}\in\overline{\mathbb{P}}$ determines $\alpha$ completely.

\begin{de}\label{defdominate}
Let $(L,\leq)$ be a linear order and let $a=\langle a_m:m\in L\rangle$, $b=\langle b_m:m\in L\rangle$ be two sequences in $\mathbb{N}_\infty$. We say that $a$ dominates $b$ if, for every $l\in L$:
\begin{equation}\label{eqdom}
\sum_{m\geq l}a_m\geq\sum_{m\geq l}b_m.
\end{equation}

For $\alpha,\beta\in{\cal A}$ we define: $\alpha\preceq\beta$ if $\beta\upharpoonright{\cal P}$ dominates $\alpha\upharpoonright{\cal P}$ for every ${\cal P}\in\overline{\mathbb{P}}$. If $\alpha\preceq\beta$ and $\beta\preceq\alpha$, we write $\alpha\approx\beta$.
\end{de}

\begin{lm}\label{pomdominate}
Let $\alpha,\beta\in{\cal A}$ be $\cal U$-closed and let ${\cal P}\in\overline{\mathbb{P}}$. The following conditions are equivalent:

(i) $\beta\upharpoonright{\cal P}$ dominates $\alpha\upharpoonright{\cal P}$;

(ii) there is a one-to-one function $g_{\cal P}:\bigcup_{u\in{\cal E}_{\cal P}}(\{u\}\times\alpha({\cal P}^u))\rightarrow\bigcup_{v\in{\cal E}_{\cal P}}(\{v\}\times\beta({\cal P}^v))$ such that $v\succeq_{\cal P}u$ whenever $g_{\cal P}(u,i)=(v,j)$;

(iii) there is a function $f_{\cal P}:\bigcup_{u\in{\cal E}_{\cal P}}(\{u\}\times\alpha({\cal P}^u))\rightarrow{\cal E}_{\cal P}$ such that $f_{\cal P}(u,i)\succeq_{\cal P}u$ for every $(u,i)\in\bigcup_{u\in{\cal E}_{\cal P}}(\{u\}\times\alpha({\cal P}^u))$ and $|f_{\cal P}^{-1}[\{v\}]|\leq\beta({\cal P}^v)$ for every $v\in{\cal E}_{\cal P}$.
\end{lm}

\dokaz (i)$\Rightarrow$(ii) Assume that $\beta\upharpoonright{\cal P}$ dominates $\alpha\upharpoonright{\cal P}$. We consider two cases.

Case 1. If $\sum_{v\in{\cal E}_{\cal P}}\beta({\cal P}^v)=n\in\mathbb{N}$, then we enumerate $\bigcup_{v\in{\cal E}_{\cal P}}(\{v\}\times\beta({\cal P}^v))=\{(v_j,l_j):j<n\}$ and $\bigcup_{u\in{\cal E}_{\cal P}}(\{u\}\times\alpha({\cal P}^u))=\{(u_j,i_j):j<k\}$ in the descending order of first coordinates. Clearly, (\ref{eqdom}) implies that $k\leq n$ and $u_j\preceq_{\cal P}v_j$ for $j<k$, so $g_{\cal P}(u_j,i_j)=(v_j,l_j)$ defines a function as desired.

Case 2. Let $m\in{\cal E}_{\cal P}$ be the maximal such that $\sum_{v\succeq_{\cal P}m}\beta({\cal P}^v)=\infty$. (By Theorem \ref{samocplus}, this sum, if infinite, must be equal to $\infty$; the maximal such $m$ exists by Example \ref{zatvex}(b).) If $m$ has an immediate successor $w$ in $({\cal E}_{\cal P},\prec_{\cal P})$, then we can define $g(u,i)\in\bigcup_{v\succeq_{\cal P}w}(\{v\}\times\beta({\cal P}^v))$ for $(u,i)\in\bigcup_{u\succeq_{\cal P}w}(\{u\}\times\alpha({\cal P}^u))$ in the same way as we defined $g_{\cal P}$ in Case 1. Otherwise, fix a descending sequence $\langle v_\xi:\xi<\gamma\rangle$ in $({\cal E}_{\cal P},\prec_{\cal P})$ such that $\inf_{\xi<\gamma}v_\xi=m$ (constructing it by recursion and using Lemma \ref{supinf}(b) at limit stages). As in Case 1, by recursion on $\xi$ we define a one-to-one function $g$ mapping each $(u,i)\in\bigcup_{u\succeq_{\cal P}v_\xi}(\{u\}\times\alpha({\cal P}^u))$ to some $g(u,i)\in\bigcup_{v\succeq_{\cal P}v_\xi}(\{v\}\times\beta({\cal P}^v))$. Thus, $|\{g(u,i):u\succ_{\cal P} m\land i<\alpha({\cal P}^u)\}|\leq\aleph_0$. Now enumerate the remaining pairs: $\bigcup_{v\succeq m}(\{v\}\times\beta({\cal P}^v))\setminus\{g(u,i):u\succ_{\cal P} m\land i<\alpha({\cal P}^u)\}=\{(v_\zeta,i_\zeta):\zeta<\infty\}$. Let $h:\bigcup_{u\preceq_{\cal P} m}(\{u\}\times\alpha({\cal P}^u))\rightarrow\infty$ be any one-to-one function. Now a desired function can be defined as
$$g_{\cal P}(u,j)=\left\{\begin{array}{ll}
g(u,j) & \mbox{if }u\succ_{\cal P} m,\\
(v_{h(u,j)},i_{h(u,j)}) & \mbox{otherwise.}
\end{array}\right.$$

(ii)$\Rightarrow$(iii) is obvious.

(iii)$\Rightarrow$(i) Let $f_{\cal P}$ be a function as in (iii). For any $u\in{\cal E}_{\cal P}$ the set $\bigcup_{w\succeq_{\cal P} u}(\{w\}\times\alpha({\cal P}^w))$ is contained in $f_{\cal P}^{-1}[\{w\in{\cal E}_{\cal P}:w\succeq_{\cal P} u\}]$, so its cardinality $\sum_{w\succeq_{\cal P} u}\alpha({\cal P}^w)$ is at most $|f_{\cal P}^{-1}[\{w\in{\cal E}_{\cal P}:w\succeq_{\cal P} u\}]|\leq\sum_{w\succeq_{\cal P} u}\beta({\cal P}^w)$.\kraj

Recall that $D^{[u,u]}_x=\{(p,k):[(p,k)]_{\approx_{\cal P}}=u\land p^{k}\zvepar x\}.$

\begin{de}\label{defalpha}
For any $x=\prod_{n\leq z}p_n^{h(n)}\in\zve{\mathbb{N}}$ as in Proposition \ref{fund}, define $\alpha_x\in{\cal A}$ as follows. For each basic ${\cal P}^u\in{\cal B}$, let $\alpha_x({\cal P}^u):=|D^{[u,u]}_x|$.
\end{de}

\begin{te}\label{zatvorenost}
For every $x\in\zve{\mathbb{N}}$, the pattern $\alpha_x$ is $\cal U$-closed.
\end{te}

\dokaz Assume that ${\cal Q}^u\in{\cal B}$, $n\in\mathbb{N}$ and $\sum_{{\cal P}^u\in\overline{A}}\alpha_x({\cal P}^u)\geq n$ for every $A\in{\cal Q}^u\cap{\cal U}$. Denote $B_A=\{(q_1,q_2,\dots,q_n)\in\zve(A\cap\mathbb{P}^{exp})^n:(\forall i\neq j)q_i\neq q_j\land(\forall i)q_i\zvepar x\}$. We prove that the family $F:=\{B_A:A\in{\cal Q}^u\cap{\cal U}\}$ has the finite intersection property. This family is closed for finite intersections, so we only need to show that each $B_A$ is nonempty. Hence let $A\in{\cal Q}^u\cap{\cal U}$; $\sum_{{\cal P}^u\in\overline{A}}\alpha_x({\cal P}^u)\geq n$ implies that there are some ${\cal P}_i^{v_i}\in\overline{A}$ and $z_i\in\mu({\cal P}_i^{v_i})$ for $1\leq i\leq n$ such that $z_i\zvepar x$. Thus $B_A\neq\emptyset$ and $F$ has the finite intersection property, so by ${\goth c}^+$-saturation we get distinct $q_1,q_2,\dots,q_n\in\bigcup_{w\succeq_{\cal Q}u}\mu({\cal Q}^w)$ such that $q_i\zvepar x$, which means that $\sum_{w\succeq_{\cal Q}u}\alpha_x({\cal Q}^w)\geq n$. Thus $\alpha_x$ is $\cal U$-closed.\kraj

We will show in Corollary \ref{postoji} that (a sort of) a converse of Theorem \ref{zatvorenost} also holds: every $\cal U$-closed pattern is $\approx$-equivalent to one of the form $\alpha_x$ for some $x$.

\begin{lm}\label{pomveza}
If $x\zvemid y$, then $\alpha_x\preceq\alpha_y$.
\end{lm}

\dokaz Let ${\cal P}\in\overline{\mathbb{P}}$. According to Definition \ref{defalpha}, to every $(u,i)\in\bigcup_{u\in{\cal E}_{\cal P}}(\{u\}\times\alpha_x({\cal P}^u))$ corresponds some $(p_{u,i},k_{u,i})\in u$, such that $p_{u,i}^{k_{u,i}}\zvepar x$ and $p_{u,i}$'s are all distinct. Let $f_{\cal P}(u,i):=[(p_{u,i},{\rm exp}_{p_{u,i}}y)]_{\approx_{\cal P}}$; clearly $f_{\cal P}(u,i)\succeq_{\cal P}u$ and $|f_{\cal P}^{-1}[\{w\}]|\leq\alpha_y({\cal P}^w)$ for every $w\in{\cal E}_{\cal P}$. By Lemma \ref{pomdominate}, the function $f_{\cal P}$ witnesses that $\alpha_x\upharpoonright{\cal P}$ is dominated by $\alpha_y\upharpoonright{\cal P}$.\kraj

\begin{te}\label{samealpha}
For any ${\cal F}\in\beta\mathbb{N}$ and any two $x,y\in\mu({\cal F})$ holds $\alpha_x=\alpha_y$.
\end{te}

\dokaz It suffices to prove that, for every ${\cal P}^u\in{\cal B}$, $\alpha_y({\cal P}^u)\geq n$ implies $\alpha_x({\cal P}^u)\geq n$. For $A\in{\cal P}^u\cap{\cal C}$ let
$$X_A:=\{m\in\mathbb{N}:(\exists q_1,\dots,q_n\in A\cap\mathbb{P}^{exp})((\forall i\neq j)q_i\neq q_j\land(\forall i)q_i\parallel m)\}.$$
Then $y\in\zve X_A$ implies $x\in\zve X_A$. Thus, the family
$$F:=\{\{(q_1,q_2,\dots,q_n)\in\zve(A\cap\mathbb{P}^{exp})^n:(\forall i\neq j)q_i\neq q_j\land(\forall i)q_i\zvepar x\}:A\in{\cal P}^u\cap{\cal C}\}$$
has the f.i.p., so by ${\goth c}^+$-saturation $x$ has distinct exact divisors $q_1,q_2,\dots,q_n\in\mu({\cal P}^u)$.\kraj

Theorem \ref{samealpha} allows us to make the following definition.

\begin{de}\label{alphault}
For ${\cal F}\in\beta\mathbb{N}$ and any $x\in\mu({\cal F})$ define $\alpha_{{\cal F}}=\alpha_x$.
\end{de}

We can now restate what we proved in Lemma \ref{pomveza} as follows.

\begin{co}\label{vezaalphamid}
(a) If ${\cal F}\widemid{\cal G}$, then $\alpha_{\cal F}\preceq\alpha_{\cal G}$.

(b) If ${\cal F}=_\sim{\cal G}$, then $\alpha_{\cal F}\approx\alpha_{\cal G}$.
\end{co}

The converse implications are false. To see this, recall that by Proposition \ref{ramsey} there is ${\cal P}\in\overline{\mathbb{P}}$ for which there are $\widemid$-incomparable ultrafilters ${\cal F},{\cal G}\in\overline{L_2}$ such that $\alpha_{\cal F}=\alpha_{\cal G}=\{({\cal P},2)\}$. The point is that a pattern $\alpha_{\cal F}$ determines all the basic ultrafilters that divide ${\cal F}$ and their multiplicity, but it does not determine its $=_\sim$-equivalence class.

\begin{ex}\label{exalphault}
(a) It was already mentioned that, for $p\in\mathbb{P}$, $p^{max}=p^\omega$. Note also that, for any ${\cal F}\in\beta\mathbb{N}$, $\alpha_{\cal F}(p^n)\leq 1$ for $n\in\omega+1$, and equality holds for at most one $n$ (because ${\cal F}=p^m\cdot p^n\cdot{\cal G}$ would actually mean that $p^{m+n}\widemid{\cal F}$). In particular, $\alpha_{\cal F}(p^\omega)=1$ is equivalent to $p^n\widemid{\cal F}$ for all $n\in\mathbb{N}$. This is why the condition $\sum_{k\in\omega+1}\alpha(p^k)\leq 1$ was included in Definition \ref{defpattern}.

(b) Recall that $MAX$ is the $\widemid$-greatest class. By \cite{So4}, Lemma 4.6, ${\cal F}\in MAX$ if and only if $m\widemid{\cal F}$ for all $m\in{\mathbb{N}}$. Let us draw this conclusion from Theorem \ref{zatvorenost}. Take any ${\cal P}\in\overline{\mathbb{P}}\setminus\mathbb{P}$ and $A\in{\cal U}\cap{\cal P}^{max}$. Since $\mathbb{P}^{exp}\in{\cal P}^{max}$, $A\cap\mathbb{P}^{exp}$ is infinite. $\{p\in\mathbb{P}:(\exists n\in\mathbb{N})p^n\in A\}$ is also infinite because ${\cal P}\notin\mathbb{P}$. Since $\alpha_{\cal F}(p^\omega)=1$ whenever $p^n\widemid{\cal F}$ for all $n\in\mathbb{N}$, this means that $\sum_{p^\omega\in\overline{A}}\alpha_{\cal F}(p^\omega)=\infty$. By $\cal U$-closedness we have $\alpha_{\cal F}({\cal P}^{max})=\infty$.

Thus, $\alpha_{MAX}$ is in the $\approx$-equivalence class of $\beta:=\{(p^\omega,1):p\in\mathbb{P}\}$, and any pattern $\alpha$ is in this class if and only if it contains $\beta$.

(c) $NMAX$ is the $\widemid$-greatest class among $\mathbb{N}$-free ultrafilters (those not divisible by any $n\in\mathbb{N}$), see Section 5 of \cite{So5}. Thus, $\alpha_{NMAX}$ is the $\approx$-equivalence class of $\{({\cal P}^{max},\infty):{\cal P}\in\overline{\mathbb{P}}\setminus\mathbb{P}\}$.

(d) If $p_n\in\zve{\mathbb{P}}$ (for $n\in\mathbb{N}$) are distinct, let ${\cal P}_n=tp(p_n/\mathbb{N})$, ${\cal F}_n={\cal P}_1{\cal P}_2\dots{\cal P}_n$ and $[{\cal G}]_\sim=\lim_{n\rightarrow\omega}{\cal F}_n$. For any prime ultrafilter ${\cal Q}\in{\rm cl}_{\cal U}(\{{\cal P}_n:n\in\mathbb{N}\})$, $\alpha_{\cal G}({\cal Q})>0$ because $\alpha_{\cal G}$ is $\cal U$-closed. Thus it may happen that none of the ${\cal F}_n$'s is divisible by $\cal Q$, but their limit is.
\end{ex}

Note that the example of $MAX$ shows that Corollary \ref{vezaalphamid} can not be strengthened: ${\cal F}=_\sim{\cal G}$ does not imply $\alpha_{\cal F}=\alpha_{\cal G}$.

\section{$F_\alpha$-sets}

In this section we describe those sets from ${\cal F}\cap{\cal U}$ that are determined by basic divisors of $\cal F$. Recall that, by Lemma \ref{utop}, the open sets of the $\cal U$-topology are of the form $\overline{A^h}$, so it will suffice to consider only such sets in the following definition. Let us call the sets of this form O-sets.

\begin{de}\label{defF}
Let $\alpha\in{\cal A}$, ${\cal P}\in\overline{\mathbb{P}}$ and $A\in{\cal P}\upharpoonright\mathbb{P}$.

For $u\in{\cal E}_{\cal P}$, an $(A,{\cal P}^{u})$-set is any O-set of the form $A^h=\bigcup_{n\in\mathbb{N}}{A_n}^n\hspace{-1mm}\uparrow$ (where $h:A\rightarrow\mathbb{N}$ and $A_n=h^{-1}[\{n\}]$) such that for some/every $(p,x)\in u$ holds $p^x\in\zve A^h$.

An $(A,{\cal P}^{w})$-set for some $w\succeq_{\cal P}u$ which is not an $(A,{\cal P}^{v})$-set for any $v\prec_{\cal P}u$ will be called an $(A,{\cal P}^{\succeq u})$-set.

An $(\alpha,A,{\cal P})$-set is any finite product of $(A,{\cal P}^u)$-sets for various $u\in{\cal E}_{\cal P}$, such that for any fixed $u$, 
if $\sum_{w\succeq_{\cal P} u}\alpha({\cal P}^w)=n\in\mathbb{N}$, then there are at most $n$ $(A,{\cal P}^{\succeq u})$-sets in the product.

An $\alpha$-set is any finite product $C_1C_2\dots C_k$ of $(\alpha,A_i,{\cal P}_i)$-sets $C_i$, with $A_i\in{\cal P}_i$, ${\cal P}_i\neq{\cal P}_j$ and $A_i\cap A_j=\emptyset$ for $i\neq j$.

Finally, $F_\alpha$ is the intersection of ${\cal U}$ with the filter generated by $\alpha$-sets.
\end{de}

\begin{ex}\label{expattern}
(a) First let $\alpha=\{({\cal P}^2,1)\}$. If $A\in{\cal P}\upharpoonright\mathbb{P}$ and $A^h=\bigcup_{n\in\mathbb{N}}{A_n}^n\hspace{-0.1cm}\uparrow$ is an $(A,{\cal P}^{\geq 2})$-set, then in order for $p^2\in\zve A^h$ and $p\notin\zve A^h$ to hold (for any $p\in\mu({\cal P})$), it is necessary that $A_2\in{\cal P}$. Every such set $A^h$ is contained in ${A_2}^2\hspace{-0.1cm}\uparrow$, so $F_\alpha$ is generated by sets of the form $A^2\hspace{-1mm}\uparrow$ for $A\in{\cal P}$.

(b) In general, let us call a pattern $\alpha$ finite if $\alpha({\cal P}^u)$ is finite for all ${\cal P}^u\in{\cal B}$ and nonzero for only finitely many ${\cal P}^u$. Finite patterns are exactly what was considered in \cite{So3}. For example, if $\alpha=\{({\cal P},2),({\cal P}^2,1),({\cal Q},1)\}$, then $F_\alpha$ is generated by the family of sets of the form $(A\hspace{-1mm}\uparrow\cdot A\hspace{-1mm}\uparrow\cdot A^2\hspace{-1mm}\uparrow)\cap B\hspace{-1mm}\uparrow=(A^{(2)}A^2B)\hspace{-0.1cm}\uparrow$, where
$$A^{(2)}A^2B=\{a_1a_2a_3^2b:a_1,a_2,a_3\in A\mbox{ are distinct }\land b\in B\}$$
for some disjoint $A\in{\cal P}\upharpoonright\mathbb{P}$, $B\in{\cal Q}\upharpoonright\mathbb{P}$.

(c) If $\alpha=\{({\cal P},\infty)\}$, then $F_\alpha$ is generated by the sets $A^{(n)}\hspace{-1mm}\uparrow$ for all $A\in{\cal P}\upharpoonright\mathbb{P}$ and all $n\in\mathbb{N}$.

(d) If $\alpha=\{({\cal P}_i,1):i\in I\}$ and ${\cal P}_i\neq{\cal P}_j$ for $i\neq j$, then $F_\alpha$ is generated by the sets $(A_{i_1}A_{i_2}\dots A_{i_n})\hspace{-1mm}\uparrow$, where $\{i_1,i_2,\dots,i_n\}\subseteq I$ and $A_{i_1}\in{\cal P}_{i_1}\upharpoonright\mathbb{P},\dots,A_{i_n}\in{\cal P}_{i_n}\upharpoonright\mathbb{P}$ are disjoint.

(e) If $\alpha=\{({\cal P}^n,1):n\in\mathbb{N}\}$ for some ${\cal P}\in\overline{\mathbb{P}}\setminus\mathbb{P}$, then $F_\alpha$ is generated by the sets $(A^{k_1}A^{k_2}\dots A^{k_m})\hspace{-1mm}\uparrow$ for some $A\in{\cal P}\upharpoonright\mathbb{P}$, some $m\in\mathbb{N}$ and $k_1,k_2,\dots,k_m\in\mathbb{N}$.

(f) If $\alpha=\{({\cal P}^u,1)\}$ for some $u\in{\cal E}_{\cal P}\setminus\omega$, then $F_\alpha$ is generated by $A^h$ for $(A,{\cal P}^{u})$-sets $A^h=\bigcup_{n\in\mathbb{N}}{A_n}^n$ and $A=\bigcup_{n\in\mathbb{N}}A_n\in{\cal P}$, such that $p^x\in\zve A^h$ whenever $(p,x)\in u$. Note that, if $u=\sup_{\beta<\gamma}u_\beta$ for some $\prec_{\cal P}$-increasing sequence $\langle u_\beta:\beta<\gamma\rangle$ in ${\cal E}_{\cal P}$ then, using Proposition \ref{sublimit}, we conclude that every $B\in{\cal P}^u\cap{\cal U}$ is also in ${\cal P}^{u_\beta}\cap{\cal U}$ for some $\beta<\gamma$. Hence every $(A,{\cal P}^u)$-set is also an $(A,{\cal P}^{u_\beta})$-set for some $\beta<\gamma$.

(g) Finally, if $\alpha=\{({\cal P}^u,2)\}$ for some $u\in{\cal E}_{\cal P}\setminus\omega$, then $F_\alpha$ is generated by sets of the form $C_1C_2$ for some $(A,{\cal P}^{u})$-sets $C_1$ and $C_2$. Note that, by Transfer, for any $(p,x),(q,y)\in u$, $p^x\in\zve C_1$ and $q^y\in\zve C_2$ imply $p^xq^y\in\zve(C_1C_2)=\zve C_1\zve C_2$. ($p^xq^y$ is a typical element of any $\mu({\cal F})$ such that $\alpha_{\cal F}=\{({\cal P}^u,2)\}$.)
\end{ex}

\begin{te}\label{efovi}
For every ${\cal F}\in\beta\mathbb{N}$, $F_{\alpha_{\cal F}}\subseteq {\cal F}\cap{\cal U}$.
\end{te}

\dokaz Take any $x\in\mu({\cal F})$ and any $\alpha_x$-set $D$; since clearly $D\in{\cal U}$, we need to prove that $D\in{\cal F}$.
Let $D=\prod_{j=1}^mB_j$ be the representation of $D$ as the product of O-sets. Since all $B_i$ are $\mid$-upward closed, 
$(\forall n\in\mathbb{N})(n\in D\Leftrightarrow(\exists b_1\in B_1)\dots(\exists b_m\in B_m)(b_1,b_2,\dots,b_m\in\mathbb{P}^{exp}\mbox{ distinct }\land (\forall j\leq m)b_j\parallel n))$ so, by Transfer, $x\in\zve D$ if and only if
$$(\exists b_1\in\zve B_1)\dots(\exists b_m\in\zve B_m)(b_1,b_2,\dots,b_m\in\zve{\mathbb{P}^{exp}}\mbox{ distinct }\land(\forall j\leq m)b_j\zvepar x).$$
By the definitions of $\alpha_x$ and $F_{\alpha_x}$, this formula is true, implying that $D\in{\cal F}$.\kraj

${\cal F}\cap{\cal U}$ is not generated by $F_{\alpha_{\cal F}}$: by Proposition \ref{ramsey}, there are ${\cal P}\in\overline{\mathbb{P}}$ such that for $\alpha=\{({\cal P},2)\}$ there are $=_\sim$-nonequivalent ultrafilters containing $F_\alpha$. However, the next result shows that $F_{\alpha_{\cal F}}$ determines the $\approx$-equivalence class of $\alpha_{\cal F}$.

\begin{te}\label{vezaalpha}
For patterns $\alpha,\beta\in{\cal A}_{cl}$, the following conditions are equivalent:

(i) $\alpha\preceq\beta$;

(ii) $F_{\alpha}\subseteq F_{\beta}$.
\end{te}

\dokaz (i)$\Rightarrow$(ii) Assume $\alpha\preceq\beta$ and let us prove that every $(\alpha,A,{\cal P})$-set $C$ is also a $(\beta,A,{\cal P})$-set. Every such $C$ is a finite product of some $(A,{\cal P}^{u_j})$-sets of the form $A^{h_j}$ for some $h_j:A\rightarrow\mathbb{N}$. By Lemma \ref{pomdominate} $\alpha\preceq\beta$ implies that there is a function $f_{\cal P}$ adjoining to each such $u_j$ some $v_j\succeq_{\cal P}u_j$, so that $A^{h_j}$ is also a $(A,{\cal P}^{v_j})$-set. For every $v\in{\cal E}_{\cal P}$ there are at most $\sum_{w\succeq_{\cal P}v}\beta({\cal P}^w)$ $(A,{\cal P}^{\succeq v})$-sets in the factorization of $C$ and hence $C$ is an $(\beta,A,{\cal P})$-set.

(ii)$\Rightarrow$(i) Assume the opposite, that $\sum_{w\succeq_{\cal P} u}\alpha({\cal P}^w)>s=\sum_{w\succeq_{\cal P} u}\beta({\cal P}^w)$ for some ${\cal P}\in\overline{\mathbb{P}}$. By $\cal U$-closedness of $\beta$ and using Lemma \ref{utop}, we find $A\subseteq\mathbb{P}$ and $h:A\rightarrow\mathbb{N}$ such that $A^h\in{\cal P}^u$ and $\sum_{{\cal Q}^v\in\overline{A^h}}\beta({\cal Q}^v)=s$. Consider the $\alpha$-set $C=(A^h)^{(s+1)}$. By (ii) there are $\beta$-sets $D_1,\dots,D_m$ such that $D_1\cap\dots\cap D_m\subseteq C$. Each $D_i$ is a product containing $s_i\leq s$ $(B,{\cal P}^{\succeq u})$-sets, and without loss of generality we may assume that $B\subseteq A$. Let us fix $p_1,p_2,\dots,p_s\in\mu({\cal P})$ and $x_1,x_2,\dots,x_s$ such that $(p_i,x_i)\in u$. $\zve D_i$ contains an element $d_i$ of the form $p_1^{x_1}p_2^{x_2}\dots p_{s_i}^{x_{s_i}}y_i$, such that, for every prime factor $q$ of $y_i$ holds $q\neq p_j$ (for $j=1,2,\dots,s$) and $q^{{\rm exp}_qy_i}\notin\zve A^h$. If we denote $d'=l.c.m.\{d_1,d_2,\dots,d_m\}$, then $d'=p_1^{x_1}p_2^{x_2}\dots p_{s'}^{x_{s'}}y'$, and again $s'\leq s$ and for every prime factor $q$ of $y'$ holds $q\neq p_j$ (for $j=1,2,\dots,s$) and $q^{{\rm exp}_qy'}\notin\zve A^h$. But $d'$ belongs to $D_1\cap\dots\cap D_m$, and it can not belong to $C$ which contains only elements with more than $s$ factors from $\zve A^h$, a contradiction.\kraj

The following example, a sequel to Example \ref{zatvex}, shows that the condition of $\cal U$-closedness in the theorem above is necessary.

\begin{ex}
(a) Let ${\cal P}\in\overline{\mathbb{P}}\setminus\mathbb{P}$. Consider the patterns $\alpha=\{({\cal Q},1):{\cal Q}\in\overline{\mathbb{P}}\setminus\mathbb{P}\}$ and $\beta=(\alpha\setminus\{({\cal P},1)\})\cup\{({\cal P},2)\}$. $\alpha$ and $\beta$ are clearly not $\approx$-equivalent. On the other hand, for any $(\beta,A,{\cal P})$-set $A^{(2)}$ there are ${\cal Q}\in\overline{A}\setminus\{{\cal P}\}$ and disjoint sets $A_1\in{\cal P}$, $A_2\in{\cal Q}$ such that $A_1\cup A_2=A$, so that $A_1A_2\subseteq A^{(2)}$. In this way we can ``replace" the second copy of $\cal P$ with $\cal Q$. Thus $A^{(2)}$ is also an $\alpha$-set, so $F_\alpha=F_\beta$. However, note that $\alpha$ is not $\cal U$-closed because in every $\cal U$-neighborhood of $\cal P$ there are (infinitely many) primes ${\cal Q}\neq{\cal P}$.

(b) Likewise, $\alpha=\{({\cal P}^n,1):n\in\mathbb{N}\}$ and $\beta=\{({\cal P}^\omega,\infty)\}$ are not $\approx$-equivalent. Still, $F_\alpha$ and $F_\beta$ generate the same filter: since there are no $(A,{\cal P}^{\omega})$-sets which are not $(A,{\cal P}^n)$-sets for some $n\in\omega$, all the $\beta$-sets are also $\alpha$-sets. Again, the explanation is that $\alpha$ is not $\cal U$-closed.
\end{ex}

\begin{lm}\label{uvecanje}
If ${\cal F},{\cal H}\in\beta\mathbb{N}\setminus\mathbb{N}$ and ${\cal P}\in\overline{\mathbb{P}}\setminus\mathbb{P}$ are such that ${\cal H}\in{\cal P}^u$ for some $u\in{\cal E}_{\cal P}$, then $\alpha_{{\cal F}\cdot{\cal H}}({\cal P}^u)=\alpha_{\cal F}({\cal P}^u)+1$ and $\alpha_{{\cal F}\cdot{\cal H}}({\cal Q}^v)=\alpha_{\cal F}({\cal Q}^v)$ for all ${\cal Q}^v\in{\cal B}\setminus\{{\cal P}^u\}$.
\end{lm}

\dokaz Let $(x,p^a)\in\mu({\cal F})\times\mu({\cal H})$ be a tensor pair. Then $x<p^a$ and $xp^a\in\mu({\cal F}\cdot{\cal H})$. By definition, $\alpha_{{\cal F}\cdot{\cal H}}({\cal P}^u)=\alpha_{xp^a}({\cal P}^u)=\alpha_x({\cal P}^u)+1=\alpha_{\cal F}({\cal P}^u)+1$ and $\alpha_{{\cal F}\cdot{\cal H}}({\cal Q}^v)=\alpha_{\cal F}({\cal Q}^v)$ for all ${\cal Q}^v\neq{\cal P}^u$.\kraj

\begin{te}\label{overunder}
Let $\beta\in{\cal A}_{cl}$ and ${\cal F}\in\beta\mathbb{N}$.

(a) If $\alpha_{\cal F}\preceq\beta$, then there is ${\cal G}\in\beta\mathbb{N}$ such that $\alpha_{\cal G}\approx\beta$ and ${\cal F}\widemid{\cal G}$.

(b) If $\beta\preceq\alpha_{\cal F}$, then there is ${\cal G}\in\beta\mathbb{N}$ such that $\alpha_{\cal G}\approx\beta$ and ${\cal G}\widemid{\cal F}$.
\end{te}

\dokaz (a) We obtain the desired ultrafilter as a limit of a $\widemid$-increasing sequence $\langle{\cal G}_\gamma:\gamma<\epsilon\rangle$. By recursion we construct this sequence, along with the sequence of respective patterns $\alpha_\delta=\alpha_{{\cal G}_\delta}$. We start with ${\cal G}_0={\cal F}$ and $\alpha_0=\alpha_{\cal F}$. Assume that $\alpha_\delta$ and ${\cal G}_\delta$ have been constructed for $\delta<\gamma$ so that $\alpha_\delta\preceq\beta$ and ${\cal F}\widemid{\cal G}_\delta$.

First we consider the successor case $\gamma=\delta+1$. If $\alpha_\delta\approx\beta$, put $\epsilon=\gamma$ and we are done. Otherwise, let ${\cal P}$ be such that $\alpha_\delta\upharpoonright{\cal P}$ does not dominate $\beta\upharpoonright{\cal P}$. We consider two cases.

Case 1. If there is $u\in{\cal E}_{\cal P}$ such that $\sum_{w\succeq_{\cal P}v}\alpha_\delta({\cal P}^w)<\sum_{w\succeq_{\cal P}v}\beta({\cal P}^w)$ for all $v\preceq_{\cal P}u$ such that $\sum_{w\succeq_{\cal P}v}\beta({\cal P}^w)<\infty$, and in particular $\sum_{w\succeq_{\cal P}u}\alpha_\delta({\cal P}^w)<\infty$, we put ${\cal G}_{\delta+1}:={\cal G}_\delta\cdot{\cal H}$ for some ${\cal H}\in{\cal P}^u$. By Lemma \ref{uvecanje} we have $\alpha_\delta\prec\alpha_{\delta+1}\preceq\beta$.

Case 2. Otherwise, there are $u,v'\in{\cal E}_{\cal P}$ such that $v'\prec_{\cal P}u$, $\sum_{w\succeq_{\cal P}v'}\alpha_\delta({\cal P}^w)=\sum_{w\succeq_{\cal P}v'}\beta({\cal P}^w)<\infty$ and $\sum_{w\succeq_{\cal P}u}\alpha_\delta({\cal P}^w)<\sum_{w\succeq_{\cal P}u}\beta({\cal P}^w)$. Define $E:=\{v\prec u:\sum_{w\succeq_{\cal P}v}\alpha_\delta({\cal P}^w)=\sum_{w\succeq_{\cal P}v}\beta({\cal P}^w)\}$ and $v_0:=\sup E$. Since there is an element $v'\in E$ such that $\sum_{w\succeq_{\cal P}v'}\alpha_\delta({\cal P}^w)$ and $\sum_{w\succeq_{\cal P}v'}\beta({\cal P}^w)$ are finite, these sums change values at finitely many $v\succeq_{\cal P}v'$, so we have $v_0\in E$ and hence $v_0\prec_{\cal P}u$. Thus we obtain $\sum_{w\succeq_{\cal P}v}\alpha_\delta({\cal P}^w)<\sum_{w\succeq_{\cal P}v}\beta({\cal P}^w)$ whenever $v_0\prec_{\cal P}v\preceq_{\cal P}u$. Now take $x\in\mu({\cal G}_\delta)$ and $q\in\mu({\cal P}^{v_0})$ such that $q\zvepar x$, choose $q'\in\mu({\cal P}^u)$ such that $q\zvemid q'$ and let $y=\frac{q'}qx$. Then ${\cal G}_{\delta+1}:=tp(y/\mathbb{N})$ is such that ${\cal G}_\delta\widemid{\cal G}_{\delta+1}$, $\alpha_{\delta+1}({\cal P}^u)=\alpha_\delta({\cal P}^u)+1$ and $\alpha_{\delta+1}({\cal P}^{v_0})=\alpha_\delta({\cal P}^{v_0})-1$, so again we have $\alpha_\delta\prec\alpha_{\delta+1}\preceq\beta$.

Finally, if $\gamma$ is a limit ordinal, let $[{\cal G}_\gamma]_\sim:=\lim_{\delta\rightarrow\gamma}{\cal G}_\delta$ and $\alpha_\gamma=\alpha_{{\cal G}_\gamma}$. Let us show that (the $\approx$-equivalence class of) $\alpha_\gamma$ is the supremum of the sequence $\langle\alpha_\delta:\delta<\gamma\rangle$ in $({\cal A}_{cl},\prec)$. Assume the opposite: that there is $\alpha'\in{\cal A}_{cl}$ such that $\alpha_\delta\preceq\alpha'$ for all $\delta<\gamma$, but $\alpha_\gamma\not\preceq\alpha'$. By Theorem \ref{vezaalpha}, there is some $A\in F_{\alpha_\gamma}\setminus F_{\alpha'}$. By Proposition \ref{sublimit} there is $\delta<\gamma$ such that $A\in F_{\alpha_\delta}$, which implies $A\in F_{\alpha'}$, a contradiction. In particular, we get $\alpha_\gamma\preceq\beta$. This concludes the construction.

(b) is proved in a similar way, with the successor case requiring only the analogue of the construction from Case 2.\kraj

It may seem that, if ${\cal F}\widemid{\cal H}$ and $\beta\in{\cal A}_{cl}$ is such that $\alpha_{\cal F}\preceq\beta\preceq\alpha_{\cal H}$, then there is ${\cal G}\in\beta\mathbb{N}$ such that $\alpha_{\cal G}\approx\beta$, ${\cal F}\widemid{\cal G}$ and ${\cal G}\widemid{\cal H}$. However, this is false, as the next example shows.

\begin{ex}
Let ${\cal P}\in\overline{\mathbb{P}}\setminus\mathbb{P}$ be arbitrary, and let $p,q\in\mu({\cal P})$ be such that $(p,q)$ is a tensor pair. Denote ${\cal F}={\cal P}^2\cdot{\cal P}$ and ${\cal H}={\cal P}^3\cdot{\cal P}^5$. Since $(p^2,q)$ and $(p^3,q^5)$ are also tensor pairs (by Proposition \ref{tensor}), it follows that $tp(p^2q/\mathbb{N})={\cal F}$ and $tp(p^3q^5/\mathbb{N})={\cal H}$.

We have $\alpha_{\cal F}=\{({\cal P},1),({\cal P}^2,1)\}$ and $\alpha_{\cal H}=\{({\cal P}^3,1),({\cal P}^5,1)\}$. Now, if $\beta=\{(({\cal P},1),({\cal P}^4,1)\}$, then clearly $\alpha_{\cal F}\preceq\beta\preceq\alpha_{\cal H}$, but there can be no ultrafilter ${\cal G}$ such that $\alpha_{\cal G}\approx\beta$ (which is in this case equivalent to $\alpha_{\cal G}=\beta$), ${\cal F}\widemid{\cal G}$ and ${\cal G}\widemid{\cal H}$. Namely, for such an ultrafilter by Proposition \ref{ekviv} there would be some $y\in\mu({\cal G})$ such that $p^2q\zvemid y$, and the only possibility is $y=p^4q$. In turn, there would be some $z\in\mu({\cal H})$ such that $p^4q\zvemid z$, so $z=p^5q^3$. However, $tp(p^5q^3/\mathbb{N})\neq tp(p^3q^5/\mathbb{N})={\cal H}$ because for $A=\{a^3b^5:a,b\in\mathbb{P}\land a<b\}$ we have $p^3q^5\in\zve A$ and $p^5q^3\notin\zve A$.
\end{ex}

As a special case of the previous theorem, we have the following.

\begin{co}\label{postoji}
For every $\cal U$-closed pattern $\beta$ there is an ultrafilter ${\cal G}$ such that $\alpha_{\cal G}\approx\beta$.
\end{co}

One may also be tempted to think that, if $\beta$ is a $\cal U$-closed pattern which is not finite (as defined in Example \ref{expattern}(b)), then there are in fact $2^{\goth c}$ ultrafilters ${\cal G}$ as described in the Corollary \ref{postoji}, because of the possibility to choose, at limit stages of the construction, different ultrafilters ${\cal W}$ on $\gamma$ and get different ${\cal G}_\gamma=\lim_{\delta\rightarrow{\cal W}}{\cal G}_\delta$. However, we saw in Example \ref{cut} that there are a basic class ${\cal P}^u$ and $p\in\mu({\cal P})$ for which $u_p$ is a singleton, say $u_p=\{x\}$. It follows that for $\beta=\{({\cal P}^u,1)\}$ there can be only one ultrafilter ${\cal F}=tp(p^x/\mathbb{N})$ such that $\alpha_{\cal F}\approx\beta$. By Theorem \ref{overunder} this ultrafilter must then be divisible by all ${\cal G}\in{\cal P}^\omega$.

\section{Self-divisible ultrafilters}

As we have seen, the pattern of an ultrafilter does not determine its $=_\sim$-divisibility class. In particular, there are many such classes it can not distinguish between. However, some properties connected to $\widemid$-divisibility can be characterized by patterns. In this section we give one such application.

\begin{de}
For ${\cal F}\in\beta\mathbb{N}$ let $D({\cal F})=\{n\in\mathbb{N}:n\widemid{\cal F}\}$.

${\cal F}$ is self-divisible if $D({\cal F})\in{\cal F}$.

${\cal F}$ is division-linear if it contains a $\mid$-chain.
\end{de}

The following notions were introduced in \cite{DLM} in order to resolve some questions left open in \cite{So6} about possible extensions of the congruence relation to ultrafilters. However, it turned out that for the self-divisibility there are many more equivalent conditions of algebraic, number-theoretic and topological nature, making this kind of ultrafilters objects of interest in their own right. For example, these nonstandard characterizations were obtained in Theorem 3.10 and Proposition 4.5 of \cite{DLM}.

\begin{pp}\label{ekvivSD}
(a) ${\cal F}$ is self-divisible if and only if for every $a,b\in\mu({\cal F})$ there is $c\in\mu({\cal F})$ such that $c\zvemid gcd(a,b)$.

(b) ${\cal F}$ is division-linear if and only if for every $a,b\in\mu({\cal F})$ holds $a\zvemid b$ or $b\zvemid a$.
\end{pp}

Our goal is to apply ideas from this paper to get another characterization of self-divisible ultrafilters. More precisely, we want to add more details to conclusions of Proposition 4.17 from \cite{DLM}. For ${\cal F}\in\beta\mathbb{N}$ denote:
\begin{eqnarray*}
{\rm exp}_p{\cal F} &=& \max\{k\in\mathbb{N}:p^k\widemid{\cal F}\} \mbox{ (for }p\in\mathbb{P}\mbox{), if it exists}\\
I_{\cal F} &=& \{p\in\mathbb{P}:p^\omega\widemid{\cal F}\}\\
J_{\cal F} &=& \{p\in\mathbb{P}:(\exists k\in\mathbb{N})k={\rm exp}_p{\cal F}\}\\
K_{\cal F} &=& \{p\in\mathbb{P}:p\nwidemid{\cal F}\}.
\end{eqnarray*}

The definition of $u_{\cal P}$ in the following theorem does not depend on the choice of $q$ by Lemma \ref{selfdivpom}.

\begin{te}\label{selfdivekv}
Let ${\cal F}\in\beta\mathbb{N}$. Define $h_0:J_{\cal F}\rightarrow\mathbb{N}$ by $h_0(p)={\rm exp}_p{\cal F}$. For any $q\in\zve J_{\cal F}$ and ${\cal P}=tp(q/\mathbb{N})$ let $u_{\cal P}\in{\cal E}_{\cal P}$ be such that $q^{\zve h_0(q)}\in\mu({\cal P}^{u_{\cal P}})$. Then:
\begin{itemize}
\item[(a)] $\alpha_{\cal F}({\cal P}^{max})=\infty$ for all ${\cal P}\in\overline{I_{\cal F}}\setminus I_{\cal F}$.

\item[(b)] $\cal F$ is self divisible if and only if the following holds:

(i) $\alpha_{\cal F}({\cal P}^{u_{\cal P}})=\infty$ and $\alpha_{\cal F}({\cal P}^w)=0$ for $w\succ_{\cal P}u_{\cal P}$, for all ${\cal P}\in\overline{J_{\cal F}}\setminus J_{\cal F}$;

(ii) $\alpha_{\cal F}({\cal P}^w)=0$ for all ${\cal P}\in\overline{K_{\cal F}}\setminus K_{\cal F}$ and all $w\in{\cal E}_{\cal P}$.
\end{itemize}
\end{te}

\dokaz (a) Let ${\cal P}\in\overline{I_{\cal F}}\setminus I_{\cal F}$; such $\cal P$ exists only if $I_{\cal F}$ is infinite. Take any $\cal U$-neighborhood $\overline{A^h}$ of ${\cal P}^{max}$. 
$A$ is infinite and belongs to $\cal P$, so for any of infinitely many $p\in I_{\cal F}\cap A$ holds $p^\omega\in\overline{A^h}$ and $\alpha_{\cal F}(p^\omega)=1$. By $\cal U$-closedness of $\alpha_{\cal F}$, $\alpha_{\cal F}({\cal P}^{max})=\infty$.

(b) The set of divisors of $\cal F$ from $\mathbb{N}$ is
\begin{eqnarray*}
D({\cal F}) &=& \{p_1^{k_1}p_2^{k_2}\dots p_m^{k_m}q_1^{l_1}q_2^{l_2}\dots q_n^{l_n}:\\
 && p_1,p_2,\dots,p_m\in A_{\cal F}\land q_1,q_2,\dots,q_n\in B_{\cal F}\land l_i\leq{\rm exp}_{q_i}{\cal F}\mbox{ for }1\leq i\leq n\}.
\end{eqnarray*}
Define $h^+:J_{\cal F}\cup K_{\cal F}\rightarrow\mathbb{N}$ by $h^+(p)=h_0(p)+1$ for $p\in J_{\cal F}$ and $h^+(p)=1$ for $p\in K_{\cal F}$. $(J_{\cal F}\cup K_{\cal F})^{h^+}$ is the complement of $D({\cal F})$, so $\cal F$ is self-divisible if and only if $(J_{\cal F}\cup K_{\cal F})^{h^+}\notin{\cal F}$. Assume that this is true.

First let ${\cal P}\in\overline{J_{\cal F}}\setminus J_{\cal F}$. Whenever $w\succ_{\cal P}u_{\cal P}$, for any generator $q^a$ of ${\cal P}^w$ we have $a\geq\zve h^+(q)$, hence $(J_{\cal F}\cup K_{\cal F})^{h^+}\in{\cal P}^w\cap{\cal U}$. 
If $\alpha_{\cal F}({\cal P}^w)>0$, $(J_{\cal F}\cup K_{\cal F})^{h^+}$ would be a $F_{\alpha_{\cal F}}$-set not in ${\cal F}$, a contradiction with Theorem \ref{efovi}.

On the other hand, $J_{\cal F}^{h_0}\cap\mathbb{P}^{exp}=\{p^a:p\in J_{\cal F}\land a\geq h_0(p)\}\in{\cal P}^{u_{\cal P}}$. Every neighborhood $\overline{A^h}$ of ${\cal P}^{u_{\cal P}}$ intersects $J_{\cal F}^{h_0}\cap\mathbb{P}^{exp}$ in infinitely many elements. Using $\cal U$-closedness as in (a), we get $\sum_{w\succeq_{\cal P}u_{\cal P}}\alpha_{\cal F}({\cal P}^w)=\infty$. Since we already concluded that $\alpha_{\cal F}({\cal P}^w)=0$ for $w\succ_{\cal P}u_{\cal P}$, it follows that $\alpha_{\cal F}({\cal P}^{u_{\cal P}})=\infty$.

Next, if there were some prime ultrafilter ${\cal P}\in\overline{K_{\cal F}}\setminus K_{\cal F}$ such that ${\cal P}\widemid{\cal F}$, then $K_{\cal F}\uparrow\in{\cal F}$. However, $D({\cal F})$ is disjoint from $K_{\cal F}\uparrow$, so it can not belong to $\cal F$.

Finally, assume that $\cal F$ satisfies conditions (i) and (ii). For any $x\in\mu({\cal F})$, Transfer implies that $x\in\zve{(J_{\cal F}\cup K_{\cal F})^{h^+}}$ if and only if $(\exists p\in\zve{K_{\cal F}})p\zvemid x$ or $(\exists p\in\zve{J_{\cal F}})p^{\zve h_0(p)+1}\zvemid x$. The first formula is false because of (iii), and the second because of (ii): remember that, by Lemma \ref{selfdivpom}, $p^{\zve h_0(p)+1}$ belongs to a different basic class than $p^{\zve h_0(p)}$. Hence $(J_{\cal F}\cup K_{\cal F})^{h^+}\notin{\cal F}$, so $D({\cal F})\in{\cal F}$.\kraj

Let us compare this result with Proposition 4.17 from \cite{DLM}. Case (1) there covers the possibility $I_{\cal F}=\emptyset$ and $J_{\cal F}$ is finite, in which $\cal F$ is self-divisible only if it is principal (by (ii)). Case (2) occurs when $J_{\cal F}$ and $K_{\cal F}$ are both finite, so every such $\cal F$ is self-divisible. Finally, Case (3) covers the rest, and this is where Theorem \ref{selfdivekv} brings new information. Note that it is easy to show that, for every partition of $\mathbb{N}$ into sets $I_{\cal F},J_{\cal F}$ and $K_{\cal F}$ and every choice of ${\rm exp}_p{\cal F}$ for $p\in J_{\cal F}$ there is an ultrafilter satisfying the conditions of the theorem.\\

Unfortunatelly, division-linearity of $\cal F$ does not depend only on the pattern $\alpha_{\cal F}$, as the following example shows.

\begin{ex}
Choose any $p,q\in\mathbb{P}$, say $p<q$. Take any $a,b\in\zve{\mathbb{N}}\setminus\mathbb{N}$ with $a<b$ such that $p^a$ and $p^b$ generate the same ultrafilter, say ${\cal G}=tp(p^a/\mathbb{N})=tp(p^b/\mathbb{N})$. Now let $c\in\zve{\mathbb{N}}$ be such that $(p^b,q^c)$ is a tensor pair (equivalently, that $(b,c)$ is a tensor pair), and let ${\cal H}=tp(q^c/\mathbb{N})$. Finally, let $d\in\mu({\cal H})$ be such that $d>c$ and $(a,d)$ is a tensor pair. Then $tp(p^aq^d/\mathbb{N})=tp(p^bq^c/\mathbb{N})={\cal G}\cdot{\cal H}$. However, $p^aq^d$ and $p^bq^c$ are $\zvemid$-incomparable since $a<b<c<d$, so ${\cal G}\cdot{\cal H}$ is not division-linear by Proposition \ref{ekvivSD}(b).

On the other hand, let us show that there is a division-linear ultrafilter with the same pattern $\alpha=\{(p^\omega,1),(q^\omega,1)\}$. Let us denote the family of all $F_\alpha$-sets (for all patterns $\alpha$) by ${\cal U}'$. $L=\{p^kq^k:k\in\mathbb{N}\}$ is a $\mid$-chain. By Theorem \ref{vezaalpha} it suffices to show that the family $S:=F_\alpha\cup\{A^c:A\in{\cal U}'\setminus F_\alpha\}\cup\{L\}$ has the finite intersection property. (In general, finding an ultrafilter with the same $F_\alpha$-sets would only show that the pattern of the obtained ultrafilter is $\approx$-equivalent to $\alpha$, but for $\alpha$ with finite support it actually means that its pattern equals $\alpha$.) A typical intersection of sets from $F_\alpha$ is of the form $(p^mq^n)\uparrow$ for some $m,n\in\mathbb{N}$. Any set $A\in{\cal U}'\setminus F_\alpha$ must not have any of the elements $p^lq^l$ for $l\in\mathbb{N}$, since otherwise it would contain $(p^lq^l)\uparrow$. Thus, for any finite subfamily of $S$, taking $k:=\max\{m,n\}$ we get an element $p^kq^k$ in the intersection of all the sets.
\end{ex}

\section{Closing remarks and open questions}

We assumed throughout the text that all sets of the form $\{1,2,\dots,z\}$ for $z\in\zve{\mathbb{N}}\setminus\mathbb{N}$ have the same cardinality; let us now consider what was gained by it. This condition implies that many other internal sets have the same cardinality. This allowed us to conclude that there is only one possible infinite value of the number of divisors of a fixed ultrafilter $\cal F$ from a given ${\cal E}_{\cal P}$-class, which we denoted by $\infty$. There are several advantages of this. First of all, this simplifies working with patterns, and does not require any additional set-theoretic assumptions. Second, it makes possible (or at least simplifies) the proofs of several results, such as Theorem \ref{samocplus} and Lemma \ref{pomdominate} that depends on it. Most notably, Theorem \ref{samealpha} is quite important: we need to know that what we learn about the pattern of $\cal F$ from some $x\in\mu({\cal F})$ is really a property of $\cal F$ and not of the particular $x$.

Here are a few questions that remain unanswered.

\begin{qu}\label{q1}
(a) Is it possible, for a given ${\cal P}\in\overline{\mathbb{P}}\setminus\mathbb{P}$, to describe precisely the order $({\cal E}_{\cal P},\prec_{\cal P})$?

(b) In particular, is $({\cal E}_{\cal P},\preceq_{\cal P})$ isomorphic to $({\cal E}_{\cal Q},\preceq_{\cal Q})$ for all nonprincipal ${\cal P}$ and ${\cal Q}$?
\end{qu}

\begin{qu}\label{q2}
Is it possible to improve Theorem \ref{overunder} (or at least Corollary \ref{postoji}) to get an ultrafilter $\cal G$ such that $\alpha_{\cal G}=\beta$?
\end{qu}



A previous version of this paper was developed under a stronger condition for which an assumption on cardinal arithmetic was needed. The author is grateful to Mauro Di Nasso for pointing out the principle $\Delta_1$, which made possible the elimination of this additional assumption. 

Funding: The author gratefully acknowledges financial support of the Science Fund of the Republic of Serbia (call PROMIS, project CLOUDS, grant no.\ 6062228 and call IDEJE, project Set-theoretic, model-theoretic and Ramsey-theoretic phenomena in mathematical structures: similarity and diversity -- SMART, grant no.\ 7750027) and of the Ministry of Science, Technological Development and Innovation of the Republic of Serbia (grant no.\ 451-03-47/2023-01/200125).

\end{document}